\documentclass[10pt,a4paper]{article}

\usepackage[utf8]{inputenc}
\usepackage[T1]{fontenc} 
\usepackage[left=2cm,right=2cm,top=2cm,bottom=2cm]{geometry}

\usepackage{authblk}

\usepackage{graphicx}
\usepackage{multirow}
\usepackage{amsmath,amssymb,amsfonts}
\usepackage{amsthm}
\usepackage{mathrsfs}
\usepackage[page]{appendix}
\usepackage{textcomp}
\usepackage{manyfoot}
\usepackage[rightcaption]{sidecap}
\usepackage{physics} 
\usepackage{lineno}
\usepackage{hyperref}
\usepackage{tabularx}
\usepackage{booktabs}

\usepackage{todonotes}

\renewcommand{\epsilon}{\varepsilon}
\renewcommand{\theta}{\vartheta}
\renewcommand{\rho}{\varrho}
\renewcommand{\phi}{\varphi}

\title{General multi-steps variable-coefficient formulation for computing quasi-periodic solutions with multiple base frequencies}

\author[1]{Junqing Wu}
\author[1]{Ling Hong}
\author[2]{Mingwu Li\thanks{limw@sustech.edu.cn}}
\author[1]{Jun Jiang\thanks{jun.jiang@mail.xjtu.edu.cn}}

\affil[1]{State Key Laboratory for Strength and Vibration, Xi’an Jiaotong University, Xi’an 710049, China}
\affil[2]{Department of Mechanics and Aerospace Engineering, Southern University of Science and Technology, Shenzhen 518055, China}

\date{\today} 

\begin{document}

\maketitle

\begin{abstract}
    Quasi-periodic solutions with multiple base frequencies exhibit the feature of $2\pi$-periodicity with respect to each of the hyper-time variables. However, it remains a challenge work, due to the lack of effective solution methods, to solve and track the quasi-periodic solutions with multiple base frequencies until now. In this work, a \textit{multi-steps variable-coefficient formulation} (\textit{m}-VCF) is proposed, which provides a unified framework to enable either harmonic balance method (HB) or collocation method (CO) or finite difference method (FD) to solve quasi-periodic solutions with multiple base frequencies. For this purpose, a method of alternating U and S domain (AUS) is also developed to efficiently evaluate the nonlinear force terms. Furthermore, a new robust phase condition is presented for all of the three methods to make them track the quasi-periodic solutions with prior unknown multiple base frequencies, while the stability of the quasi-periodic solutions is assessed by mean of Lyapunov exponents. The feasibility of the constructed methods under the above framework is verified by application to three nonlinear systems.
\end{abstract}

\vspace{0.5cm} 


\section{Introduction}\label{sec:introduction}

As one of the common steady-state motions of nonlinear dynamical systems, quasi-periodic solutions have been reported in many literatures of various fields \cite{JW1,JW2,JW3,JW4,JW5,JW6}. In past decades, a variety of numerical methods have been proposed and continuously under development. Most of them are suitable for the quasiperiodic solutions only with two base frequencies and although few of them claim to be applicable to those with multiple base frequencies, they are often complicate in operation without detailed introduction of the procedures and application examples. 

As known, the numerical methods for the quasi-periodic solutions can be roughly classified into three categories according to how the solutions are expressed in time domain. For \textit{the first category of methods}, the quasi-periodic solution is expressed in the physical time domain. The most familiar and used method is the time integration one (TI), which is time consuming but not reliable in obtaining complete quasi-periodic solutions due to the non-repeatability of quasi-periodic solutions. For a special case when all the base frequencies are known and they have an approximate greatest common divisor, the quasi-periodic solution can be then approximately considered as a periodic solution with a single approximate frequency. Thus, the harmonic balance method (HB) \cite{JW7,JW8} or the time variational method (TVM) \cite{JW9} can be used to solve it. Furthermore, a higher-order Poincaré map, as an extension of the shooting method for periodic solutions, is proposed based on the boundary condition of quasi-periodic solutions in \cite{JW10,JW11,JW12}. However, due to the non-repeatability of quasi-periodic solutions, the mapping of an initial point of a quasi-periodic solution can only be approximated under interpolation through a long run of time integration. As known, the multiple harmonic balance method (MHB) \cite{JW13} is initially used to approximate the quasi-periodic solution by the truncated Fourier series in time domain, where frequency components are defined as the linear combinations of the incommensurable base frequencies. However, the nonlinear force terms cannot be reliably evaluated through a simple extension of the classical alternating frequency-time (AFT) method due to the non-repeatability of nonlinear force terms in the time domain \cite{JW14,JW15}.

For \textit{the second category of methods}, the quasi-periodic solution is expressed in the hyper-time domain, where the hyper-time variables are the product of the base frequencies and physical time \cite{JW16,JW17,JW18}. Given that the quasi-periodic solution exhibits $2\pi$-periodicity with respect to each hyper-time variable $\tau_i=\omega_it$, the solution can be uniquely determined within a hyper-cubic domain. Only in this hyper-time domain can MHB be promoted and applied because the nonlinear force terms can be now effectively evaluated by a new AFT \cite{JW19,JW20,JW21}. According to the definition of MHB, it can compute the quasi-periodic solutions with multiple base frequencies. However, the detailed selection of harmonic order is not actually provided in \cite{JW19,JW22,JW23} and the evaluation of the nonlinear force terms by new AFT is time-consuming because it involves the operations of the large-scale matrices, such as those of multi-dimensional discrete Fourier transforms (\textit{m}-DFT) and their inverse transforms (\textit{m}-\textit{i}DFT). It is known that the incremental harmonic balance (IHB) method with two-time scales \cite{JW24,JW25} is a variant of MHB for the two base frequencies. On the other hand, the variable-coefficient harmonic balance method (VCHB) \cite{JW20,JW26} and twice harmonic balance method (THB) \cite{JW27} are another idea for solving the quasi-periodic solutions that are approximated by a set of periodic functions with $\tau_1$ along with the corresponding variable-coefficients, which is function with respect to $\tau_2$. Then a second harmonic balance procedure is used to solve the variable-coefficients. In \cite{JW18,JW28}, the finite difference method (FD) has been also extended to compute the quasi-periodic solutions with multiple base frequencies. However, only examples with two base frequencies were presented. Also, a spectral collocation method \cite{29} was proposed for the delay differential systems to solve the quasi-periodic solutions with two base frequencies.

For \textit{the third category of methods}, the quasi-periodic solution is expressed in a variant of the hyper-time domain by defining some new variables $\varphi_1=\tau_1$ and $\varphi_i=\tau_i+\omega_i/\omega_1\tau_1,\,i\neq1$. The $2\pi$-periodicity with respect to the hyper-time variable $\tau_1$ is thereby transformed into a new boundary condition with respect to the variable $\varphi_1$, which is restricted between 0 and $2\pi$. As can be seen, the quasi-periodic solutions having  -periodicity with respect to $\varphi_i,\,i\neq1$ are also preserved. For instance, the method proposed by Dankowicz and Schilder \cite{JW30} uses the Fourier series expansion along with a set of variable-coefficients to approximate the quasi-periodic solutions, where the Fourier series are the $2\pi$-periodic functions with respect to $\varphi_2$ and the variable-coefficients are the functions with respect to $\varphi_1$. Then the multiple trajectories constructed by the collocation method coupled through the boundary conditions based on Fourier series expansion are used to compute the quasi-periodic solutions with two base frequencies. Recently, a general purpose toolbox based on this method has been developed \cite{JW31}. In this aspect, a shooting method in the hyper-time domain is also devised recently in \cite{JW28} by the boundary condition with respect to the variable $\varphi_1$. The residual of this shooting method is approximated by the cubic splines and the examples of quasi-periodic solutions with only two base frequencies are presented. 

From the above discussion it is found that there is still lack on the general framework to solve the quasi-periodic solutions with multiple base frequencies on basis of the traditional methods for the periodic solutions. Therefore, a \textit{multi-steps variable-coefficient formulation} (\textit{m}-VCF) is proposed in this work with the aim to present a unified framework for calculating the quasi-periodic solutions with multiple base frequencies that are applicable to the three methods, namely, the frequency-domain method: harmonic balance method (HB) \cite{JW32,JW33,JW34,JW35}, and the time-domain methods: finite difference method (FD) \cite{JW18,JW36} and collocation method (CO) \cite{JW37,JW38}. It is well known that the three methods are very good at solving periodic solutions. So, they are adopted to approximate the quasi-periodic solutions with multiple base frequencies in the hyper-time domain. In the implementation, the solutions are firstly approximated by a set of $2\pi$-periodic functions in terms of the first hyper-time variable $\tau_1$ along with the corresponding variable-coefficients that are function of the rest hyper-time variables. The periodic functions are then discretized, resulting in a set of ordinary differential equations, whose solutions are just the variable-coefficients. Repeating the above operations multiple times, a set of nonlinear algebraic equations for a set of constant coefficients will be eventually obtained and the zero roots of the equations provide the quasi-periodic solution with multiple base frequencies. Under the unified formulation, quasi-periodic solutions are eventually discretized through the combinations of three methods in a mixed domain (called U domain in \textit{m}-VCF).

It is emphasized here that under the present framework each set of the periodic functions with respect to a $\tau_i$ can be independently solved by anyone of the above mentioned three numerical methods for periodic solutions. For example, with the \textit{m}-VCF, a quasi-periodic solution with two base frequencies can be solved through HB+HB or HB+FD or HB+CO etc. The rest combinations are given as FD+HB, FD+FD, FD+CO, CO+HB, CO+FD and CO+CO. To facilitate the evaluation on the nonlinear force terms in the hyper-time domain, the Alternating Frequency-Time approach \cite{JW15,JW20} is extended in the form of an alternating U and S domain method (AUS). As known, it is challenge to track the full branches of the quasi-periodic solutions by varying system parameters in order to explore their stability and bifurcation. For this purpose, a simple and robust phase condition proposed in the authors’ recent work \cite{JW39} are now developed to track the quasi-periodic solution with prior unknown multiple base frequencies applicable for all the three methods. Moreover, the general-purpose toolbox of the m-VCF is incorporated into the MATLAB program of continuation in NLvib \cite{JW32}.

The rest of this paper is organized as follows. The \textit{m}-VCF is proposed in Section \ref{sec2}, which sets up a unified framework for constructing periodic functions with variable-coefficients and their discretization, and can integrate with the three methods for periodic solutions. Then, in Section \ref{sec3}, the AUS method is used to evaluate nonlinear force terms, which is based on two defined concepts of U domain and S domain. In Section \ref{sec4}, a robust phase condition is further developed for quasi-periodic solutions with multiple prior unknown base frequencies. An initialization technique of quasi-periodic solutions with two base frequencies is also briefly reviewed in this section. In Section \ref{sec5}, the feasibility of the proposed methods is demonstrated by applying to three nonlinear systems. Finally, conclusions are drawn in Section \ref{sec6}.
\section{The multi-steps variable-coefficient formulation}\label{sec2}

In this work, a second-order nonlinear differential system with $n$-DOFs is considered:
\begin{equation}
	\begin{aligned}
		&\mathbf{R}^{0}\left(\mathbf{Z}^{0}\left(t\right),t\right)=\mathbf{M}^0\ddot{\mathbf{Z}}^{0}\left(t\right)+\mathbf{D}^{0}\dot{\mathbf{Z}}^{0}\left(t\right)+\mathbf{K}^{0}\mathbf{Z}^{0}\left(t\right)+\boldsymbol{\Theta}^{0}\left[\mathbf{F}^{0}(\mathbf{Z}^{0},\dot{\mathbf{Z}}^{0})-\mathbf{E}^{0}\left(\boldsymbol{\omega}_1^et\right)\right]=\mathbf{0}
	\end{aligned},\label{eq1}
\end{equation}
where over-dots denote the derivatives with respect to physical time $t\in[0,+\infty)$; $\mathbf{M}^0,\,\mathbf{D}^0,\,\mathbf{K}^0$ are the generalized mass, the viscous damping and the stiffness matrices with the size of $n\times n$, respectively; $\boldsymbol{\Theta}^{0}$ is the force distribution matrix generally being $n\times n$ identity matrix; $\mathbf{F}^{0}(\mathbf{Z}^{0},\dot{\mathbf{Z}}^{0})\in\mathbb{R}^{n\times1}$ represents the vector of nonlinear forces that depend on generalized displacements $\mathbf{Z}^{0}$ and generalized velocities $\dot{\mathbf{Z}}^{0}$; $\mathbf{E}^{0}(\boldsymbol{\omega}_1^et)$ is the vector of quasi-periodic excitation with frequencies of $\boldsymbol{\omega}_1^e:=\left[\omega_1,\omega_2,\cdots,\omega_e\right]^\mathrm{T}$, which has the  -periodicity with respect to each argument $\omega_it$, i.e., $\mathbf{E}^0(\cdots,\omega_it,\cdots)=\mathbf{E}^0(\cdots,\omega_it+2\pi,\cdots),\,i=1,\cdots,e$. A special case is $e=0$, when the system is not under the external excitation, i.e., when $\mathbf{E}^0\equiv\mathbf{0}$.

Assuming that the system exhibits the quasi-periodic solutions with $d$ base frequencies, i.e., $\mathbf{Z}^0(\boldsymbol{\omega}_1^dt):\mathbb{R}^{d\times1}\times[0,+\infty)\mapsto\mathbb{R}^{n\times1},\,d\geq e$. Here, $\boldsymbol{\omega}_1^d$ is $\left[\omega_1,\omega_2,\cdots,\omega_d\right]^\mathrm{T}$ according to the definition of $\boldsymbol{\omega}_1^e$. It should be emphasized that the frequencies $\omega_{i},\,e<i\leq d$ of the quasi-periodic solutions $\mathbf{Z}^0(\boldsymbol{\omega}_1^dt)$ is not known as a priori. Thus, the setting in Eq. \eqref{eq1} encompasses the cases of autonomous, periodically-forced, and quasi-periodically-forced motions. By introducing the hyper-time variables $\boldsymbol{\tau}_1^d=\boldsymbol{\omega}_1^dt\in\mathbb{S}^d$, the quasi-periodic solution is transformed into $\mathbf{Z}^0(\boldsymbol{\tau}_1^d){:}\mathbb{S}^d\to\mathbb{R}^{n\times1}$ and it can be uniquely determined in a hyper-cubic subspace $\mathbb{T}^d=[0,2\pi)^d\subset\mathbb{S}^d$ of the hyper-time domain.  

To introduce multi-steps variable-coefficient formulation (\textit{m}-VCF), consider a set of nonlinear differential equations, which exhibit the quasi-periodic solutions $\mathbf{Z}^{i-1}(\boldsymbol{\tau}_i^d):\mathbb{S}^{d-i+1}\to\mathbb{R}^{n\mathrm{U}_1^{i-1}\times1}$ in the hyper-time domain:
\begin{equation}
	\begin{aligned}
		\mathbf{R}^{i-1}(\mathbf{Z}^{i-1}(\boldsymbol{\tau}_i^d),\boldsymbol{\tau}_i^d)&=\mathbf{M}^{i-1}\ddot{\mathbf{Z}}^{i-1}(\boldsymbol{\tau}_i^d)+\mathbf{D}^{i-1}\dot{\mathbf{Z}}^{i-1}(\boldsymbol{\tau}_i^d)+\mathbf{K}^{i-1}\mathbf{Z}^{i-1}(\boldsymbol{\tau}_i^d)\\&+\boldsymbol{\Theta}^{i-1}\times[\mathbf{F}^{i-1}(\mathbf{Z}^{i-1},\dot{\mathbf{Z}}^{i-1})-\mathbf{E}^{i-1}(\boldsymbol{\tau}_i^e)]=\mathbf{0}
	\end{aligned},\label{eq2}
\end{equation}
where $\mathbf{M}^{i-1},\,\mathbf{D}^{i-1},\,\mathbf{K}^{i-1},\,\boldsymbol{\Theta}^{i-1}\in\mathbb{R}^{n\mathrm{U}_{1}^{i-1}\times n\mathrm{U}_{1}^{i-1}}$ are the generalized mass, the viscous damping, the stiffness and the force distribution matrices, respectively; And over-dots denote the derivatives with respect to physical time $t\in[0,+\infty)$; $\mathbf{F}^{i-1},\mathbf{E}^{i-1}\in\mathbb{R}^{n\mathrm{U}_{1}^{i-1}\times1}$ are the vectors of nonlinear force terms and excitation.

In the \textit{i}-th step of the \textit{m}-VCF, the quasi-periodic solution $\mathbf{Z}^{i-1}(\boldsymbol{\tau}_i^d)$ and its derivatives are approximated by choosing any one of the numerical methods for periodic solutions (HB, CO and FD):
\begin{equation}
	\begin{aligned}
		&\mathbf{Z}^{i-1}(\boldsymbol{\tau}_{i}^{d})\approx[\mathbf{I}_{n\mathrm{U}_{1}^{i-1}}\otimes\boldsymbol{\Phi}_{i}(k_{i},\tau_{i})]\mathbf{Z}^{i}(k_{i},\boldsymbol{\tau}_{i+1}^{d})
		\\&\dot{\mathbf{Z}}^{i-1}(\tau_{i}^{d})\approx[\mathbf{I}_{n\mathrm{U}_{1}^{i-1}}\otimes\omega\boldsymbol{\Phi}^{\prime}_{i}(k_{i},\tau_{i})]\mathbf{Z}^{i}(k_{i},\boldsymbol{\tau}_{i+1}^{d})+[\mathbf{I}_{n\mathrm{U}_1^{i-1}}\otimes\boldsymbol{\Phi}_i(k_i,\tau_i)]\dot{\mathbf{Z}}^i(k_i,\boldsymbol{\tau}_{i+1}^d)
		\\&\ddot{\mathbf{Z}}^{i-1}(\boldsymbol{\tau}_{i}^{d})\approx[\mathbf{I}_{n\mathrm{U}_{1}^{i-1}}\otimes\omega_{i}{}^{2}\boldsymbol{\Phi}^{\prime\prime}_{i}(k_{i},\tau_{i})]\mathbf{Z}^{i}(k_{i},\boldsymbol{\tau}_{i+1}^{d})+[\mathbf{I}_{n\mathrm{U}_1^{i-1}}\otimes2\omega_i\boldsymbol{\Phi}^{\prime}_{i}(k_i,\tau_i)]\dot{\mathbf{Z}}^i(k_i,\boldsymbol{\tau}_{i+1}^d)\\&+[\mathbf{I}_{n\mathrm{U}_1^{i-1}}\otimes\boldsymbol{\Phi}_i(k_i,\tau_i)]\ddot{\mathbf{Z}}^i(k_i,\boldsymbol{\tau}_{i+1}^d)
	\end{aligned},\label{eq3}
\end{equation}
where $\otimes$ is the Kronecker tensor product; $\mathbf{I}_{n\mathrm{U}_{1}^{i-1}}$ is $n\mathrm{U}_{1}^{i-1}\times\mathrm{U}_{1}^{i-1}$ identity matrix; $\boldsymbol{\Phi}_i(k_i,\tau_i){:}\mathbb{U}_i\times\mathbb{T}\to\mathbb{R}^{1\times\mathrm{U}_i},\,i=1,2,\cdots,d$ collects $2\pi$-periodic functions with respect to $\tau_i$. Here, $k_i$ is called as \textit{k}-parameter to suppose these periodic solutions, where the U domain $\mathbb{U}_i$ has $\mathrm{U}$ elements. $(\cdot)^{\prime}$ and $(\cdot)^{\prime\prime}$ represent the derivatives with respect to $\tau_i$. And $\mathbf{Z}^i(k_i,\boldsymbol{\tau}_{i+1}^d):\mathbb{U}_i\times\mathbb{T}^{d-i}\to\mathbb{R}^{n\mathrm{U}_1^i\times1}$ are the corresponding variable-coefficients with respect to $k_i$ and $\boldsymbol{\tau}_{i+1}^d$.

Substituting Eq. \eqref{eq3} into Eq. \eqref{eq2} and discretizing the periodic functions $\boldsymbol{\Phi}_i,\,\boldsymbol{\Phi}_i^{\prime},\,\boldsymbol{\Phi}_i^{\prime\prime}$ yields a new set of nonlinear differential equations:
\begin{equation}
	\begin{aligned}
		&\mathbf{R}^{i}(\mathbf{Z}^{i}(\boldsymbol{\tau}_{i+1}^d),\boldsymbol{\tau}_{i+1}^d)=\mathbf{M}^{i}\ddot{\mathbf{Z}}^{i}(\boldsymbol{\tau}_{i+1}^d)+\mathbf{D}^{i}\dot{\mathbf{Z}}^{i}(\boldsymbol{\tau}_{i+1}^d)+\mathbf{K}^{i}\mathbf{Z}^{i}(\boldsymbol{\tau}_{i+1}^d)+\boldsymbol{\Theta}^{i}\times[\mathbf{F}^{i}(\mathbf{Z}^{i},\dot{\mathbf{Z}}^{i})-\mathbf{E}^{i}(\boldsymbol{\tau}_{i+1}^e)]=\mathbf{0}
	\end{aligned},\label{eq4}
\end{equation}
with
\begin{equation}
	\begin{aligned}
		&\mathbf{M}^{i}=\mathbf{M}^{i-1}\otimes\boldsymbol{\Upsilon}_{i}{}^{0}
		\\&\mathbf{D}^{i}=2\mathbf{M}^{i-1}\otimes\omega_{i}\boldsymbol{\Upsilon}_{i}{}^{1}+\mathbf{D}^{i-1}\otimes\boldsymbol{\Upsilon}_{i}{}^{0}
		\\&\mathbf{K}^{i}=\mathbf{M}^{i-1}\otimes\omega_{i}^{2}\boldsymbol{\Upsilon}_{i}{}^{2}+\mathbf{D}^{i-1}\otimes\omega_{i}\boldsymbol{\Upsilon}_{i}{}^{1}+\mathbf{K}^{i-1}\otimes\boldsymbol{\Upsilon}_{i}{}^{0}
		\\&\boldsymbol{\Theta}^{i}=\boldsymbol{\Theta}^{i-1}\otimes\boldsymbol{\Upsilon}_{i}{}^{0}
	\end{aligned},\label{eq5}
\end{equation}
where $\boldsymbol{\Upsilon}_i{}^0(k_i),\,\boldsymbol{\Upsilon}_i{}^1(k_i),\,\boldsymbol{\Upsilon}_i{}^2(k_i):\mathbb{U}_i\to\mathbb{R}^{\mathrm{U}_i\times\mathrm{U}_i}$ are three constant matrices because $k_i$ needs to be chosen before computing solutions. Considering the previous $i-1$ steps of the \textit{m}-VCF, Eq. \eqref{eq4} should be in the following form:
\begin{equation}
	\begin{aligned}
		\mathbf{R}^{i}(\mathbf{Z}^{i}(\boldsymbol{k}_{1}^i,\boldsymbol{\tau}_{i+1}^d),\boldsymbol{\omega}_{1}^i,\boldsymbol{\tau}_{i+1}^d)&=\mathbf{M}^{i}(\boldsymbol{\omega}_{1}^i)\ddot{\mathbf{Z}}^{i}(\boldsymbol{k}_{1}^i,\boldsymbol{\tau}_{i+1}^d)
		+\mathbf{D}^{i}(\boldsymbol{\omega}_{1}^i)\dot{\mathbf{Z}}^{i}(\boldsymbol{k}_{1}^i,\boldsymbol{\tau}_{i+1}^d)+\mathbf{K}^{i}(\boldsymbol{\omega}_{1}^i)\mathbf{Z}^{i}(\boldsymbol{k}_{1}^i,\boldsymbol{\tau}_{i+1}^d)
		\\&+\boldsymbol{\Theta}^{i}(\boldsymbol{k}_{1}^i)\times[\mathbf{F}^{i}(\mathbf{Z}^{i},\dot{\mathbf{Z}}^{i})-\mathbf{E}^{i}(\boldsymbol{k}_{1}^i,\boldsymbol{\tau}_{i+1}^e)]=\mathbf{0}
	\end{aligned},\label{eq6}
\end{equation}
where $\boldsymbol{k}_{1}^i$ collects the k-parameter $\boldsymbol{k}_{1}^i=\left[k_1,k_2,\cdots,k_i\right]^\mathrm{T}$. And the notation $\mathrm{U}_{1}^{i-1}$ in Eq. \eqref{eq3} is equal to $\prod_{j=1}^{i-1}\mathrm{U}_j$. After $d$ steps of the \textit{m}-VCF, there is a set of nonlinear algebraic equations $\mathbf{R}^d(\mathbf{Z}^d(\boldsymbol{k}_1^d),\boldsymbol{\omega}_1^d)=\mathbf{0}\in\mathbb{R}^{n\mathrm{U}_1^d\times1}$ with respect to a set of constant coefficients $\mathbf{Z}^d(\boldsymbol{k}_1^d):\mathbb{U}_1^d\to\mathbb{R}^{n\mathrm{U}_1^d\times1}$ and base frequencies $\boldsymbol{\omega}_1^d\in\mathbb{R}^{d\times1}$:
\begin{equation}
	\begin{aligned}
		&\mathbf{R}^d(\mathbf{Z}^d(\boldsymbol{k}_1^d),\boldsymbol{\omega}_1^d)=\mathbf{K}^d(\boldsymbol{\omega}_1^d)\mathbf{Z}^d(\boldsymbol{k}_1^d)
		+\boldsymbol{\Theta}^d(\boldsymbol{k}_1^d)\times[\mathbf{F}^d(\mathbf{Z}^d)-\mathbf{E}^d(\boldsymbol{k}_1^d)]=\mathbf{0}
	\end{aligned},\label{eq7}
\end{equation}
Form the above procedure, it can be found that the quasi-periodic solution $\mathbf{Z}^{0}(\boldsymbol{\tau}_{1}^{d}){:}\mathbb{S}^{d}\to\mathbb{R}^{n\times1}$ in the hyper-time domain is eventually transformed into a set of coefficients $\mathbf{Z}^{d}(\boldsymbol{k}_{1}^{d}){:}\mathbb{U}_1^{d}\to\mathbb{R}^{n\mathrm{U}_i^d\times1}$ in the U domain by the \textit{m}-VCF. And the zero roots of $\mathbf{R}^d(\mathbf{Z}^d(\boldsymbol{k}_1^d),\boldsymbol{\omega}_1^d)=\mathbf{0}$ can return back to quasi-periodic solutions:
\begin{equation}
	\begin{aligned}
		\mathbf{Z}^0(t)&\approx\mathbf{Z}^0(\boldsymbol{\tau}_1^d)
		=[\mathbf{I}_n\otimes\boldsymbol{\Phi}_1(k_1,\tau_1)\otimes\cdots\otimes\boldsymbol{\Phi}_d(k_d,\tau_d)]\mathbf{Z}^d(\boldsymbol{k}_1^d)
	\end{aligned},\label{eq8}
\end{equation}
The solution curves along the continuation parameter $p$ are calculated by the continuation technique, i.e., the tangent prediction and orthogonal corrections: 
\begin{equation}
	\begin{aligned}
		&\begin{bmatrix}\partial_\mathbf{y}\mathbf{R}^d|_{(\mathbf{y},p)_j}&\partial_p\mathbf{R}^d|_{(\mathbf{y},p)_j}\\\Delta\mathbf{y}_{j-1}^\mathrm{T}&\Delta{p}_{j-1}\end{bmatrix}
		\begin{bmatrix}\Delta\mathbf{y}_j\\\Delta{p}_j\end{bmatrix}=\begin{bmatrix}\mathbf{0}\\1\end{bmatrix}
		\\&\mathbf{y}_{j+1}^0=\mathbf{y}_j+s_j\Delta\mathbf{y}_j,\,{p}_{j+1}^0={p}_j+s_j\Delta{p}_j
	\end{aligned},\label{eq9}
\end{equation}
\begin{equation} 
	\begin{aligned}
		&\begin{bmatrix}\partial_\mathbf{y}\mathbf{R}^d|_{(\mathbf{y},p)_{j+1}^k}&\partial_p\mathbf{R}^d|_{(\mathbf{y},p)_{j+1}^k}\\\Delta\mathbf{y}_j^\mathrm{T}&\Delta{p}_j\end{bmatrix}
		\begin{bmatrix}\delta\mathbf{y}_{j+1}^{k+1}\\\delta{p}_{j+1}^{k+1}\end{bmatrix}=\begin{bmatrix}-\mathbf{R}^d|_{(\mathbf{y},p)_{j+1}^k}\\0\end{bmatrix}
		\\&\mathbf{y}_{j+1}^{k+1}=\mathbf{y}_{j+1}^k+\delta\mathbf{y}_{j+1}^{k+1},\,{p}_{j+1}^{k+1}={p}_{j+1}^k+\delta{p}_{j+1}^{k+1}
	\end{aligned},\label{eq10}
\end{equation}
where $\partial_\mathbf{y}\mathbf{R}^d,\partial_p\mathbf{R}^d$ stands for the derivative of $\mathbf{R}^d$ with respect to $\mathbf{y},\,p$. Here, $\mathbf{y}$ is defined as $\mathbf{y}=[\mathbf{Z}^d;\boldsymbol{\omega}_1^d]$. $[\Delta\mathbf{y}_j;\Delta p_j]=[\Delta\mathbf{y}_j;\Delta p_j]/\left|\Delta p_j\right|$ is the tangent vector $\mathbf{t}_{j}$, $s_j$ is the step of continuation. $(\mathbf{y},p)_{j+1}^{k}$ is considered as a solution of Eq. \eqref{eq10} until $\left\|\mathbf{R}^d\right\|_2<\varepsilon$ is satisfied, where $\varepsilon$ is a user-defined accuracy. 

Here, $\partial_{{\mathbf{Z}^{d}}}\mathbf{R}^{d}=\mathbf{K}^{d}\left(\boldsymbol{\omega}_{1}^{d}\right)+\boldsymbol{\Theta}^{d}\partial_{{\mathbf{Z}^{d}}}\mathbf{F}^{d}$ and $\partial_{\omega_i}\mathbf{R}^d=\partial_{\omega_i}\mathbf{K}^d\mathbf{Z}^d$ with 
\begin{equation}
	\begin{aligned}
		&\partial_{{\omega_{i}}}\mathbf{K}^{d}=2\mathbf{M}^{i-1}\otimes\boldsymbol{\Upsilon}_{i}{}^{1}\otimes\sum_{k=i+1}^{d} \omega_k \left(\otimes_{l=i+1}^{d}\boldsymbol{\Upsilon}_{l}{}^{\lceil0,1\rceil}\right)+
		\left(\mathbf{M}^{i-1}\otimes2\omega_{i}\boldsymbol{\Upsilon}_{i}{}^{2}+\mathbf{D}^{i-1}\otimes\boldsymbol{\Upsilon}_{i}{}^{1}\right)\otimes\left(\otimes_{i=i+1}^{d}\boldsymbol{\Upsilon}_{i}{}^{0}\right)
	\end{aligned},\label{eq11}
\end{equation}
where $\otimes_{l=i+1}^{d}\boldsymbol{\Upsilon}_{l}{}^{0}=\boldsymbol{\Upsilon}_{i+1}{}^{0}\otimes...\otimes\boldsymbol{\Upsilon}_{d}{}^{0}$ and $\boldsymbol{\Upsilon}_l^{\lceil0,1\rceil}=\boldsymbol{\Upsilon}_l{}^1$ when $k=l$ , otherwise $\boldsymbol{\Upsilon}_l{}^{\lceil0,1\rceil}=\boldsymbol{\Upsilon}_l{}^0$. And the derivative $\partial_p\mathbf{R}^d$ can be determined by the numerical difference method if $p$ is not $\omega_1$. 

Moreover, the details of the construction of functions $\boldsymbol{\Phi}_i(k_i,\tau_i)$ and their discretization based on HB, CO and FD are introduced in Appendix \ref{appA}, \ref{appB} and \ref{appC}, respectively.

\section{The Alternating U and S domain method for nonlinear terms}\label{sec3}

Like the quasi-periodic solutions in Eq. \eqref{eq8}, the nonlinear force terms can be also expressed by:
\begin{equation}
	\begin{aligned}
		&\mathbf{F}^0(\mathbf{Z}^0(t),\dot{\mathbf{Z}}^0(t))\approx\mathbf{F}^0(\mathbf{Z}^0(\boldsymbol{\tau}_1^d),\dot{\mathbf{Z}}^0(\boldsymbol{\tau}_1^d))
		=[\mathbf{I}_n\otimes\boldsymbol{\Phi}_1(k_1,\tau_1)\otimes\cdots\otimes\boldsymbol{\Phi}_d(k_d,\tau_d)]\mathbf{F}^d(\mathbf{Z}^d(\boldsymbol{k}_1^d))
	\end{aligned},\label{eq12}
\end{equation}
However, the coefficients $\mathbf{F}^d(\mathbf{Z}^d)$ and their derivative $\partial_{{\mathbf{Z}^{d}}}\mathbf{F}^{d}$ are usually hard to directly determine in the U domain. To address this, the alternating U and S domain (AUS) is proposed to indirectly evaluate them based on the S domain $\overline{\boldsymbol{\tau}}_{1}^{d}\in\overline{\mathbb{S}}^{d}\subset\mathbb{S}^{d}$, which is the discretized subspace of the hyper-time domain:
\begin{equation}
	\begin{aligned}
		&\overline{\boldsymbol{\tau}}_1^i=\left[\boldsymbol{O}_{\mathrm{S}_i}\otimes\overline{\boldsymbol{\tau}}_1^{i-1},\overline{\boldsymbol{\tau}}_i\otimes\boldsymbol{O}_{\mathrm{S}_1^{i-1}}\right]\in\mathbb{R}^{\mathrm{S}_1^i\times i}
		\\&\overline{\boldsymbol{\tau}}_i^d=\left[\boldsymbol{O}_{\mathrm{S}_{i+1}^d}\otimes\overline{\boldsymbol{\tau}}_i,\overline{\boldsymbol{\tau}}_{i+1}^d\otimes\boldsymbol{O}_{\mathrm{S}_i}\right]\in\mathbb{R}^{\mathrm{S}_i^d\times(d-i+1)}
	\end{aligned},\label{eq13}	
\end{equation}
where $\overline{\boldsymbol{\tau}}_i=:\left[\tau_{i;1},\cdots,\tau_{i;\mathrm{S}_i}\right]^\mathrm{T}\subset[0,2\pi)$ with $\tau_{i;j}=2\left(j-1\right)\pi/\mathrm{S}_i,\,j=1,2,\cdots,\mathrm{S}_{i}$, $\mathrm{S}_{i}$ is the number of equal-spaced points with respect to $\tau_i$, $\boldsymbol{O}_{s_1^i},\boldsymbol{O}_{s_{i+1}^d}, \boldsymbol{O}_{s_i}$ are all-ones column vectors in size of $\mathrm{S}_{1}^{i}=\prod_{k=1}^{i}\mathrm{S}_{k}$, $\mathrm{S}_{i+1}^{d}=\prod_{k=d}^{i+1}\mathrm{S}_{k}$ and $\mathrm{S}_{i}$, respectively. Noted that $\mathrm{S}_{i}$ is set as $\mathrm{U}_{i}$ when the function $\boldsymbol{\Phi}_i(k_i,\tau_i)$ is constructed by CO or FD. For HB, $\mathrm{S}_{i}\geq\mathrm{U}_{i}$ is used to avoid aliasing due to mixture of the high-order into the low-order harmonics.

The ASU method is divided into two processes, namely the process from U domain to S domain (UTS):
\begin{equation}
	\begin{aligned}	
		&\mathbb{Z}^d(\boldsymbol{k}_1^d)\Rightarrow\begin{bmatrix}\mathbf{Z}^{d-1}(\boldsymbol{k}_1^{d-1},\overline{\boldsymbol{\tau}}_d)\\\dot{\mathbf{Z}}^{d-1}(\boldsymbol{k}_1^{d-1},\overline{\boldsymbol{\tau}}_d)\end{bmatrix}\Rightarrow\cdots\Rightarrow\begin{bmatrix}\mathbf{Z}^{i}(\boldsymbol{k}_1^{i},\overline{\boldsymbol{\tau}}_{i+1}^d)\\\dot{\mathbf{Z}}^{i}(\boldsymbol{k}_1^{i},\overline{\boldsymbol{\tau}}_{i+1}^d)\end{bmatrix}
		\overset{i-\mathrm{UTS}}{\operatorname*{\Longrightarrow}}\begin{bmatrix}\mathbf{Z}^{i-1}(\boldsymbol{k}_1^{i-1},\overline{\boldsymbol{\tau}}_i^d)\\\dot{\mathbf{Z}}^{i-1}(\boldsymbol{k}_1^{i-1},\overline{\boldsymbol{\tau}}_i^d)\end{bmatrix}\Rightarrow\cdots\Rightarrow\begin{bmatrix}\mathbf{Z}^{0}(\overline{\boldsymbol{\tau}}_1^d)\\\dot{\mathbf{Z}}^{0}(\overline{\boldsymbol{\tau}}_1^d)\end{bmatrix}
	\end{aligned},\label{eq14}	
\end{equation}
and the process from S domain to U domain (STU) containing two parts:
\begin{equation}
	\begin{aligned}	
		&\mathbb{F}^0(\overline{\boldsymbol{\tau}}_1^d)\Rightarrow\mathbb{F}^1(k_1,\overline{\boldsymbol{\tau}}_2^d)\Rightarrow\cdots\Rightarrow\mathbb{F}^{i-1}(\boldsymbol{k}_1^{i-1},\overline{\boldsymbol{\tau}}_i^d)
		\overset{i-\mathrm{STU-}1}{\operatorname*{\Longrightarrow}}\mathbb{F}^{i}(\boldsymbol{k}_1^{i},\overline{\boldsymbol{\tau}}_{i+1}^d)\Rightarrow\cdots\Rightarrow\mathbb{F}^{d}(\boldsymbol{k}_1^{d})
	\end{aligned},\label{eq15}	
\end{equation}
\begin{equation}
	\begin{aligned}	
		&\begin{bmatrix}\partial_{\mathbf{Z}^0}\mathbf{F}^0\\\partial_{\dot{\mathbf{Z}}^0}\mathbf{F}^0\end{bmatrix}\Rightarrow\begin{bmatrix}\partial_{\mathbf{Z}^1}\mathbf{F}^1\\\partial_{\dot{\mathbf{Z}}^1}\mathbf{F}^1\end{bmatrix}\Rightarrow\cdots\Rightarrow\begin{bmatrix}\partial_{\mathbf{Z}^{i-1}}\mathbf{F}^{i-1}\\\partial_{\dot{\mathbf{Z}}^{i-1}}\mathbf{F}^{i-1}\end{bmatrix}
		\overset{i-\mathrm{STU-}2}{\operatorname*{\Longrightarrow}}\begin{bmatrix}\partial_{\mathbf{Z}^{i}}\mathbf{F}^{i}\\\partial_{\dot{\mathbf{Z}}^{i}}\mathbf{F}^{i}\end{bmatrix}\Rightarrow\cdots\Rightarrow\begin{bmatrix}\partial_{\mathbf{Z}^{d-1}}\mathbf{F}^{d-1}\\\partial_{\dot{\mathbf{Z}}^{d-1}}\mathbf{F}^{d-1}\end{bmatrix}\Rightarrow\partial_{\mathbf{Z}^{d}}\mathbf{F}^{d}
	\end{aligned},\label{eq16}	
\end{equation}
The following are the details of these two processes.

\subsection{The process from U domain to S domain}\label{sec3.1}

The process of UTS is to transform constant-coefficients $\mathbf{Z}^d(\boldsymbol{k}_1^d)$ in the U domain $\boldsymbol{k}_1^d\in\mathbb{U}_1^d$ into the discrete displacements $\mathbf{Z}^{0}(\overline{\boldsymbol{\tau}}_1^d)$ and velocities $\dot{\mathbf{Z}}^{0}(\overline{\boldsymbol{\tau}}_1^d)$ in the S domain $\overline{\boldsymbol{\tau}}_{1}^{d}\in\overline{\mathbb{S}}^{d}$. Here is the \textit{i}-UTS step in Eq. \eqref{eq14}:
\begin{equation}
	\underbrace{\begin{bmatrix}\mathbf{Z}^i(\boldsymbol{k}_1^i,\boldsymbol{\overline{\tau}}_{i+1}^d)\\\dot{\mathbf{Z}}^i(\boldsymbol{k}_1^i,\boldsymbol{\overline{\tau}}_{i+1}^d)\end{bmatrix}}_{\mathbb{R}^{n\mathrm{S}_{i+1}^d\times\mathrm{U}_1^i}}\overset{i-\mathrm{UTS}}{\operatorname*{\Longrightarrow}}\underbrace{\begin{bmatrix}\mathbf{Z}^{i-1}(\boldsymbol{k}_1^{i-1},\boldsymbol{\overline{\tau}}_i^d)\\\dot{\mathbf{Z}}^{i-1}(\boldsymbol{k}_1^{i-1},\boldsymbol{\overline{\tau}}_i^d)\end{bmatrix}}_{\mathbb{R}^{n\mathrm{S}_i^d\times\mathrm{U}_1^{i-1}}},\label{eq17}
\end{equation}
According to the definition of the Eq. \eqref{eq17}, the $\mathbf{Z}^d(\boldsymbol{k}_1^d)$ in Eq. \eqref{eq14} is the matrix with size of $n\times\mathrm{U}_1^d$, which is constructed by $\mathbf{Z}^d:=[\mathbf{Z}_1{}^d,\mathbf{Z}_2{}^d,\cdotp\cdotp\cdotp,\mathbf{Z}_n{}^d]^\mathrm{T}$ with $\mathbf{Z}_l{}^d\in\mathbb{R}^{\mathrm{U}_1^d\times1},\,l=1,\cdots,n$ being the coefficients of \textit{l}-th displacement $Z_l{}^0(t)$ in the U domain. Actually, the $\mathbf{Z}^d(\boldsymbol{k}_1^d)$ in Eq. \eqref{eq7} is $\mathbf{Z}^d:=[\mathbf{Z}_1{}^d;\mathbf{Z}_2{}^d;\cdotp\cdotp\cdotp;\mathbf{Z}_n{}^d]\in\mathbb{R}^{n\mathrm{U}_1^d\times1}$. $\mathrm{S}_{d+1}^d$ and $\mathrm{U}_{1}^0$ are equal to $1$ in Eq. \eqref{eq17}.

Then introduce the procedures in \textit{i}-UTS step. Firstly, the matrices $\mathbf{Z}^i,\dot{\mathbf{Z}}^i\in\mathbb{R}^{n\mathrm{S}_{i+1}^d\times\mathrm{U}_1^i}$ need to be transformed into $\mathbf{Z}^i,\dot{\mathbf{Z}}^i\in\mathbb{R}^{\mathrm{U}_i\times n\mathrm{S}_{i+1}^d\mathrm{U}_1^{i-1}}$:
\begin{equation}
	\begin{aligned}	
		&\mathbb{R}^{{n\mathrm{S}_{i+1}^{d}\times\mathrm{U}_{1}^{i}}}\;\overset{r}{\operatorname*{\rightarrow}}\;\mathbb{R}^{{n\mathrm{S}_{i+1}^{d}\times\mathrm{U}_{i}\times\mathrm{U}_{1}^{i-1}}}\;\overset{t}{\operatorname*{\rightarrow}}\;
		\mathbb{R}^{{\mathrm{U}_{i}\times n\mathrm{S}_{i+1}^{d}\times\mathrm{U}_{1}^{i-1}}}\;\overset{r}{\operatorname*{\rightarrow}}\;\mathbb{R}^{{\mathrm{U}_{i}\times n\mathrm{S}_{i+1}^{d}\mathrm{U}_{1}^{i-1}}}
	\end{aligned},\label{eq18}	
\end{equation}
where $\overset{r}{\operatorname*{\to}},\,\overset{t}{\operatorname*{\to}}$ represent the procedure of reshape and transpose in MATLAB, respectively.

The second procedure is the transformation from the domain of $(\boldsymbol{k}_1^i,\overline{\boldsymbol{\tau}}_{i+1}^d)$ to $(\boldsymbol{k}_1^{i-1},\overline{\boldsymbol{\tau}}_{i}^d)$:
\begin{equation}
	\begin{aligned}
		& \mathbf{Z}^{i-1}=\boldsymbol{\Gamma}_{i}{}^0\mathbf{Z}^{i}\\  & \dot{\mathbf{Z}}^{i-1}=\omega_{i}\boldsymbol{\Gamma}_{i}{}^1\mathbf{Z}^{i}+\boldsymbol{\Gamma}_{i}{}^0\dot{\mathbf{Z}}^{i}
	\end{aligned},\label{eq19}
\end{equation}
where the transform matrices $\boldsymbol{\Gamma}_{i}{}^0,\,\boldsymbol{\Gamma}_{i}{}^1\in\mathbb{R}^{\mathrm{S}_i\times\mathrm{U}_i}$ for three methods (HB, CO and FD) are also introduced in Appendix \ref{appA}, \ref{appB} and \ref{appC}, respectively. 

Thirdly, a transformation is used to obtain $\mathbf{Z}^{i-1},\dot{\mathbf{Z}}^{i-1}\in\mathbb{R}^{n\mathrm{S}_i^d\times\mathrm{U}_1^{i-1}}$ by $\mathbf{Z}^{i-1},\dot{\mathbf{Z}}^{i-1}\in\mathbb{R}^{\mathrm{S}_{i}\times n\mathrm{S}_{i+1}^{d}\mathrm{U}_{1}^{i-1}}$ in Eq. \eqref{eq19}:
\begin{equation}
	\mathbb{R}^{\mathrm{S}_i\times n\mathrm{S}_{i+1}^d\mathrm{U}_1^{i-1}}\;\overset{r}{\operatorname*{\rightarrow}}\;\mathbb{R}^{n\mathrm{S}_i^d\times\mathrm{U}_1^{i-1}},\label{eq20}
\end{equation}
After \textit{d} steps of UTS, $\mathbf{Z}^0(\overline{\boldsymbol{\tau}}_1^d),\,\dot{\mathbf{Z}}^0(\overline{\boldsymbol{\tau}}_1^d)\in\mathbb{R}^{n\mathrm{S}_1^d\times1}$ are eventually obtained. The discrete nonlinear forces $\mathbf{F}^0(\overline{\boldsymbol{\tau}}_1^d)\in\mathbb{R}^{n\mathrm{S}_1^d\times1}$ are computed by $\mathbf{F}^0(\mathbf{Z}^0,\dot{\mathbf{Z}}^0)$ in Eq. \eqref{eq1}. 

\subsection{The first part process from S domain to U domain}\label{sec3.2}

The details of \textit{i}-STU-1 in Eq. \eqref{eq15} are show as follows:
\begin{equation}
	\underbrace{\mathbf{F}^{i-1}(\boldsymbol{k}_1^{i-1},\overline{\boldsymbol{\tau}}_i^d)}_{\mathbb{R}^{n\mathrm{U}_1^{i-1}\times\mathrm{S}_i^d}}\overset{i-\mathrm{STU-}1}{\operatorname*{\Longrightarrow}}\underbrace{\mathbf{F}^i(\boldsymbol{k}_1^i,\overline{\boldsymbol{\tau}}_{i+1}^d)}_{\mathbb{R}^{n\mathrm{U}_1^i\times\mathrm{S}_{i+1}^d}},\label{eq21}
\end{equation}
According to the definition of the Eq. \eqref{eq21}, the $\mathbf{F}^0(\overline{\boldsymbol{\tau}}_1^d)$ in Eq. \eqref{eq15} is the matrix with size of $n\times\mathrm{S}_1^d$, which is constructed by $\mathbf{F}^d:=[\mathbf{F}_1{}^0,\mathbf{F}_2{}^0,\cdotp\cdotp\cdotp,\mathbf{F}_n{}^0]^\mathrm{T}$ with $\mathbf{F}_j{}^0\in\mathbb{R}^{\mathrm{S}_1^d\times1},\,j=1,\cdots,n$ being \textit{j}-th nonlinear force term $F_j{}^0(t)$ in the S domain. Here, the matrix $\mathbf{F}^0(\overline{\boldsymbol{\tau}}_1^d)$ obtained in the above subsection is $\mathbf{F}^d:=[\mathbf{F}_1{}^0;\mathbf{F}_2{}^0;\cdotp\cdotp\cdotp;\mathbf{F}_n{}^0]\in\mathbb{R}^{n\mathrm{S}_1^d\times1}$. $\mathrm{S}_{d+1}^d$ and $\mathrm{U}_{1}^0$ are also equal to $1$ in Eq. \eqref{eq21}.

Firstly, the matrices $\mathbf{F}^{i-1}\in\mathbb{R}^{n\mathrm{U}_1^{i-1}\times\mathrm{S}_i^d}$ need to be transformed into $\mathbf{F}^{i-1}\in\mathbb{R}^{\mathrm{S}_i\times n\mathrm{U}_1^{i-1}\mathrm{S}_{i+1}^d}$:
\begin{equation}
	\begin{aligned}	
		&\mathbb{R}^{{n\mathrm{U}_{1}^{i-1}\times\mathrm{S}_{i}^{d}}}\;\overset{r}{\operatorname*{\rightarrow}}\;\mathbb{R}^{{n\mathrm{U}_{1}^{i-1}\times\mathrm{S}_{i}\times\mathrm{S}_{i+1}^{d}}}\;\overset{t}{\operatorname*{\rightarrow}}\;\mathbb{R}^{{\mathrm{S}_{i}\times n\mathrm{U}_{1}^{i-1}\times\mathrm{S}_{i+1}^{d}}}\;\overset{r}{\operatorname*{\rightarrow}}\;\mathbb{R}^{{{{\mathrm{S}_{i}\times n\mathrm{U}_{1}^{i-1}\mathrm{S}_{i+1}^{d}}}}}
	\end{aligned},\label{eq22}	
\end{equation}

The second procedure is the transformation from the domain of $(\boldsymbol{k}_1^{i-1},\overline{\boldsymbol{\tau}}_{i}^d)$ to $(\boldsymbol{k}_1^i,\overline{\boldsymbol{\tau}}_{i+1}^d)$:
\begin{equation}
	\begin{aligned}
		\mathbf{F}^i=\boldsymbol{\Gamma}_i{}^{0,-1}\mathbf{F}^{i-1}
	\end{aligned},\label{eq23}
\end{equation}
where the transform matrix $\boldsymbol{\Gamma}_{i}{}^{0,-1}\in\mathbb{R}^{\mathrm{U}_i\times\mathrm{S}_i}$ for three methods (HB, CO and FD) is also introduced in Appendix \ref{appA}, \ref{appB} and \ref{appC}, respectively.

Thirdly, a transformation is used to obtain $\mathbf{F}^{i}\in\mathbb{R}^{n\mathrm{U}_1^{i}\times\mathrm{S}_{i+1}^d}$ by $\mathbf{F}^i\in\mathbb{R}^{\mathrm{U}_i\times n\mathrm{U}_1^{i-1}\mathrm{S}_{i+1}^d}$ in Eq. \eqref{eq23}:
\begin{equation}
	\mathbb{R}^{\mathrm{U}_i\times n\mathrm{U}_1^{i-1}\mathrm{S}_{i+1}^d}\;\overset{r}{\operatorname*{\rightarrow}}\;\mathbb{R}^{n\mathrm{U}_1^i\times\mathrm{S}_{i+1}^d},\label{eq24}
\end{equation}
After \textit{d} steps of STU-1, $\mathbf{F}^d(\boldsymbol{k}_1^d)\in\mathbb{R}^{n\mathrm{U}_1^d\times1}$ are eventually obtained. 

\subsection{The second part process from S domain to U domain}\label{sec3.3}

The details of \textit{i}-STU-2 in Eq. \eqref{eq16} are show as follows:
\begin{equation}
	\underbrace{\begin{bmatrix}\partial_{\mathbf{Z}^{i-1}}\mathbf{F}^{i-1}(\boldsymbol{k}_1^{i-1},\overline{\boldsymbol{\tau}}_i^d)\\\partial_{\dot{\mathbf{Z}}^{i-1}}\mathbf{F}^{i-1}(\boldsymbol{k}_1^{i-1},\overline{\boldsymbol{\tau}}_i^d)\end{bmatrix}}_{\mathbb{R}^{(n\mathrm{U}_1^{i-1})^2\times\mathrm{S}_i^d}}\overset{i-\mathrm{STU}-2}{\operatorname*{\operatorname*{\Rightarrow}}}\underbrace{\begin{bmatrix}\partial_{\mathbf{Z}^i}\mathbf{F}^i(\boldsymbol{k}_1^i,\overline{\boldsymbol{\tau}}_{i+1}^d)\\\partial_{\dot{\mathbf{Z}}^i}\mathbf{F}^i(\boldsymbol{k}_1^i,\overline{\boldsymbol{\tau}}_{i+1}^d)\end{bmatrix}}_{\mathbb{R}^{(n\mathrm{U}_1^i)^2\times\mathrm{S}_i^d}},\label{eq25}
\end{equation}

According to the definition of the Eq. \eqref{eq25}, the $\partial_{\mathbf{z}^0}\mathbf{F}^0,\,\partial_{\dot{\mathbf{z}}^0}\mathbf{F}^0$ in Eq. \eqref{eq16} is the matrix with size of $n^2\times\mathrm{S}_1^d$, which is constructed by $\partial_{\mathbf{Z}^{0}}\mathbf{F}^{0}:=[\partial_{\mathbf{Z}^{0}}\mathbf{F}_{1}{}^{0},\partial_{\mathbf{Z}^{0}}\mathbf{F}_{2}{}^{0},\cdotp\cdotp\cdotp,\partial_{\mathbf{Z}^{0}}\mathbf{F}_{n}{}^{0}]^{\mathrm{T}}$ where the element is $\partial_{\mathbf{Z}^{0}}\mathbf{F}_{j}{}^{0}=[\partial_{{\mathbf{Z}}_{1}{}^{0}}\mathbf{F}_{j}{}^{0},\partial_{\mathbf{Z}_{2}{}^{0}}\mathbf{F}_{j}{}^{0},\cdotp\cdotp\cdotp,\partial_{\mathbf{Z}_{n}{}^{0}}\mathbf{F}_{j}{}^{0}]\in\mathbb{R}^{\mathrm{S}_{1}^{d}\times n},j=1,\cdotp\cdotp\cdotp,n$ with $\partial_{\mathbf{Z}_{l}{}^{0}}\mathbf{F}_{j}{}^{0}$ being the derivatives of the \textit{j}-th nonlinear force term $F_j{}^0(t)$ with respect to the \textit{l}-th displacement $Z_l{}^0(t)$ in the S domain. $\mathrm{S}_{d+1}^d$ and $\mathrm{U}_{1}^0$ are also equal to $1$ in Eq. \eqref{eq25}.

Firstly, the matrices $\mathbb{R}^{(n\mathrm{U}_1^{i-1})^2\times\mathrm{S}_i^d}$ need to be transformed into $\mathbb{R}^{\mathrm{S}_{i}\times\mathrm{S}_{i}\times(n\mathrm{U}_{1}^{i-1})^{2}\mathrm{S}_{i+1}^{d}}$:

\begin{equation}
	\begin{aligned} &\mathbb{R}^{\left(n\mathrm{U}_1^{i-1}\right)^2\times\mathrm{S}_i^d}\;\overset{r}{\operatorname*{\rightarrow}}\;\mathbb{R}^{\left(n\mathrm{U}_1^{i-1}\right)^2\times\mathrm{S}_{i}\times\mathrm{S}_{i+1}^d}\;\overset{t}{\operatorname*{\rightarrow}}\;\mathbb{R}^{\mathrm{S}_{i}\times\left(n\mathrm{U}_1^{i-1}\right)^2\times\mathrm{S}_{i+1}^d}\\&\overset{r}{\operatorname*{\rightarrow}}\;\mathbb{R}^{\mathrm{S}_{i}\times1\times\left(n\mathrm{U}_1^{i-1}\right)^2\mathrm{S}_{i+1}^d}\;\overset{diag}{\operatorname*{\rightarrow}}\;\mathbb{R}^{\mathrm{S}_{i}\times\mathrm{S}_{i}\times\left(n\mathrm{U}_1^{i-1}\right)^2\mathrm{S}_{i+1}^d}
	\end{aligned},\label{eq26}
\end{equation}

where the elements on the diagonal of each page in matrix $\mathbb{R}^{\mathrm{S}_{i}\times\mathrm{S}_{i}\times\left(n\mathrm{U}_1^{i-1}\right)^2\mathrm{S}_{i+1}^d}$ are the vector of that of $\mathbb{R}^{\mathrm{S}_{i}\times1\times\left(n\mathrm{U}_1^{i-1}\right)^2\mathrm{S}_{i+1}^d}$. The second procedure is the transformation from the domain of $(\boldsymbol{k}_1^{i-1},\overline{\boldsymbol{\tau}}_{i}^d)$ to $(\boldsymbol{k}_1^i,\overline{\boldsymbol{\tau}}_{i+1}^d)$:
\begin{equation}
	\begin{aligned}
		&\partial_{\mathbf{Z}^i}\mathbf{F}^i=\boldsymbol{\Gamma}_i{}^{0,-1}\left(\partial_{\mathbf{Z}^{i-1}}\mathbf{F}^{i-1}\boldsymbol{\Gamma}_i{}^0+\omega_i\partial_{\dot{\mathbf{Z}}^{i-1}}\mathbf{F}^{i-1}\boldsymbol{\Gamma}_i{}^1\right)
		\\&\partial_{\dot{\mathbf{Z}}^i}\mathbf{F}^i=\boldsymbol{\Gamma}_i{}^{0,-1}\left(\partial_{\dot{\mathbf{Z}}^{i-1}}\mathbf{F}^{i-1}\boldsymbol{\Gamma}_i{}^0\right)
	\end{aligned},\label{eq27}
\end{equation}

Thirdly, a transformation is used to obtain $\partial_{\mathbf{Z}^i}\mathbf{F}^i,\,\partial_{\dot{\mathbf{Z}}^i}\mathbf{F}^i\in\mathbb{R}^{\left(n\mathrm{U}_1^i\right)^2\times\mathrm{S}_{i+1}^d}$ by $\partial_{{\mathbf{Z}}^i}{\mathbf{F}}^i,\,\partial_{{\dot{\mathbf{Z}}}^i}{\mathbf{F}}^i\in\mathbb{R}^{\mathrm{U}_i\times\mathrm{U}_i\times\left(n\mathrm{U}_1^{i-1}\right)^2\mathrm{S}_{i+1}^d}$ in Eq. \eqref{eq27}:
\begin{equation}
	\begin{aligned}
		&\mathbb{R}^{\mathrm{U}_i\times\mathrm{U}_i\times\left(n\mathrm{U}_1^{i-1}\right)^2\mathrm{S}_{i+1}^d}\;\overset{r}{\operatorname*{\rightarrow}}\;\mathbb{R}^{\mathrm{U}_i\times n\mathrm{U}_1^i\times n\mathrm{U}_1^{i-1}\mathrm{S}_{i+1}^d}\;\overset{t}{\operatorname*{\rightarrow}}\;\mathbb{R}^{n\mathrm{U}_1^i\times\mathrm{U}_i\times n\mathrm{U}_1^{i-1}\mathrm{S}_{i+1}^d}\;\overset{r}{\operatorname*{\rightarrow}}\;\mathbb{R}^{\left(n\mathrm{U}_1^i\right)^2\times\mathrm{S}_{i+1}^d}
	\end{aligned},\label{eq28}
\end{equation}
After \textit{d} steps of STU-2, the derivative $\partial_{\mathbf{Z}^d}\mathbf{F}^d\in\mathbb{R}^{n\mathrm{U}_1^d\times n\mathrm{U}_1^d}$ is constructed by $\partial_{\mathbf{Z}^d}\mathbf{F}^d\in\mathbb{R}^{\left(n\mathrm{U}_1^d\right)^2\times1}$ in Eq. \eqref{eq27}:
\begin{equation}
	\mathbb{R}^{\left(n\mathrm{U}_1^d\right)^2\times1}\;\overset{r}{\operatorname*{\rightarrow}}\;\mathbb{R}^{n\mathrm{U}_1^d\times n\mathrm{U}_1^d}\;\overset{t}{\operatorname*{\rightarrow}}\;\mathbb{R}^{n\mathrm{U}_1^d\times n\mathrm{U}_1^d}
	,\label{eq29}
\end{equation}

\section{Phase condition and Initialization of the solution after NS}\label{sec4}

\subsection{Phase condition}\label{sec4.1}

If $\omega_{i}$ is not explicit in $\mathbf{E}^0\left(\boldsymbol{\omega}_1^p,t\right)$ , the phase condition is required to fixed the phase of hyper-time variable $\tau_{i}$. A robust phase condition based on the procedure of continuation in Eqs. \eqref{eq9} and \eqref{eq10} is introduced in this work which is developed by authors in \cite{JW39}: 
\begin{equation}
	\begin{aligned}	
		&PC_i=\frac1{\left(2\pi\right)^d}\int_0^{2\pi}\cdots\int_0^{2\pi}\left[\frac{\partial\mathbf{Z}^{0,\mathrm{T}}}{\partial\tau_i}\,\frac{\partial^2\mathbf{Z}^{0,\mathrm{T}}}{\partial\tau_i^2}\right]\begin{bmatrix}\mathbf{Z}^0{}_*\\\frac{\partial\mathbf{Z}^0{}_*}{\partial\tau_i}\end{bmatrix}
		\mathrm{d}\tau_1\cdots\mathrm{d}\tau_d\equiv 0
	\end{aligned},\label{eq30}
\end{equation}
where $\mathbf{Z}^0{}_*$ is believed as the a vector of quasi-periodic terms whose coefficients are $\left.\Delta\mathbf{Z}^d\right|_j$ in Eq. \eqref{eq9} or $\left.\delta\mathbf{Z}^d\right|_{j+1}^{k+1}$ in Eq. \eqref{eq10}. By using Fourier-Galerkin procedure in HB or Riemann Integral method in CT and FD, Eq. \eqref{eq30} is rewritten as:
\begin{equation}
	PC_i=\left[\left(\boldsymbol{\Lambda}_i{}^1\mathbf{Z}^d\right)^\mathrm{T}+\left(\boldsymbol{\Lambda}_i{}^2\mathbf{Z}^d\right)^\mathrm{T}\boldsymbol{\Lambda}_i{}^1\right]\mathbf{Z}^d{}_*\equiv 0,\label{eq31}
\end{equation}
where  $\mathbf{Z}^d{}_*$ is  $\left.\Delta\mathbf{Z}^d\right|_j$ in Eq. \eqref{eq9} or $\left.\delta\mathbf{Z}^d\right|_{j+1}^{k+1}$ in Eq. \eqref{eq10}. And $\boldsymbol{\Lambda}_i{}^{1}, \boldsymbol{\Lambda}_i{}^{2}\in\mathbb{R}^{n\mathrm{U}_1^d\times n\mathrm{U}_1^d}$ are $\boldsymbol{\Lambda}_i{}^{1}=\mathbf{I}_{n\mathrm{U}_{1}^{i-1}}\otimes\boldsymbol{\Upsilon}_{\tau_{i}}{}^{1}\otimes \mathbf{I}_{\mathrm{U}_{i+1}^{d}}$ and $\boldsymbol{\Lambda}_i{}^{2}=\mathbf{I}_{n\mathrm{U}_{1}^{i-1}}\otimes\boldsymbol{\Upsilon}_{\tau_{i}}{}^{2}\otimes \mathbf{I}_{\mathrm{U}_{i+1}^{d}}$, respectively. $\boldsymbol{\Upsilon}_{\tau_{i}}{}^{1},\boldsymbol{\Upsilon}_{\tau_{i}}{}^{2}\in\mathbb{R}^{\mathrm{U}_i\times\mathrm{U}_i}$ are also constant matrices, whose details are shows in Appendix \ref{appA}, \ref{appB} and \ref{appC}.

\subsection{Initialization of the quasi-periodic solution after NS}\label{sec4.2}

Neimark-Sacker bifurcation (NS) may be observed in the continuation of periodic solution, where the periodic solution becomes unstable and a quasi-periodic torus occurs at point of NSB. The stabilities of quasi-periodic solutions with multiple base frequencies are evaluated by the method of Lyapunov exponents, which is introduced in Appendix \ref{appD}. In this work, an initial guess of this tours is adapted from the derivation in \cite{JW31}. 

Consider the first order of nonlinear differential equations of Eq. \eqref{eq1}
\begin{equation}
	\dot{\mathbf{x}}=\mathbf{f}\left(t,\mathbf{x}\right),\;\mathbf{x}:=\left[\mathbf{Z}^0;\dot{\mathbf{Z}^0}\right],\label{eq32}
\end{equation}
with a periodic solution $\mathbf{x}_p\left(t\right)$ with period $T_1=\frac{2\pi}{\omega_1}$. Given a solution $\mathbf{x}\left(t\right)=\mathbf{x}_p\left(t\right)+\Delta\mathbf{x}\left(t\right)$ with $\Delta\mathbf{x}\left(t\right)$ being a perturbation, the Eq. \eqref{eq32} is transformed into: 
\begin{equation}
	\Delta\dot{\mathbf{x}}(t)\approx\mathbf{J}(t)\Delta\mathbf{x}(t),\;\mathbf{J}(t)=\frac{\partial\mathbf{f}(t,\mathbf{x})}{\partial\mathbf{x}(t)},\label{eq33}
\end{equation}
Let $\boldsymbol{\Psi}(t)\in\mathbb{R}^{2n\times2n}$ as a perturbation system, then $\boldsymbol{\dot{\Psi}}(t)\approx\mathbf{J}(t)\boldsymbol{\Psi}(t)$, $\boldsymbol{\Psi}(0)=\mathbf{I}_{2n}$. The transition matrix $\boldsymbol{M}(t,0)$ is defined as $\boldsymbol{\Psi}(t)$ and the monodromy matrix is $\boldsymbol{M}(T_1,0)=\boldsymbol{\Psi}(T_1)$. When NS occurs, the monodromy matrix $\boldsymbol{M}(T_1,0)$ has at least one pair of conjugate imaginary eigenvalues $e^{\mathrm{i}\alpha},\,e^{-\mathrm{i}\alpha}$ and their eigenvectors $\upsilon(0),\bar{\upsilon}(0)$, where $\Re\{\upsilon(0)\}$ and $\Im\{\upsilon(0)\}$ span an invariant subspace to the monodromy matrix $\boldsymbol{M}(T_1,0)$. The eigenvectors  $\upsilon(t)$ of the monodromy matrix $\boldsymbol{M}(t+T_1,t)$ is calculated by $\upsilon(t)=\mathrm{e}^{-\mathrm{i}\alpha\frac t{T_1}}\boldsymbol{M}(t,0)\upsilon(0)$, $t\in\left[0,T_1\right)$, which is a periodic function with period $T_1$. Introducing the hyper-time variable $\tau_1=\omega_1t$ yields $\upsilon(\tau_1)=\mathrm{e}^{-\mathrm{i}\alpha\frac{\tau_1}{2\pi}}\boldsymbol{M}(\tau_1,0)\upsilon(0)$, $\tau_1\in\left[0,2\pi\right)$. Then, the initial guess of quasi-periodic torus is given by: 
\begin{equation}
	\begin{aligned}	
		&\mathbf{x}_{qp}\left(\tau_1,\tau_2\right)=\mathbf{x}_p\left(\tau_1\right)+\varepsilon\left[\mathrm{cos}\left(\tau_2\right)\Re\left\{\upsilon\left(\tau_1\right)\right\}+\mathrm{sin}\left(\tau_2\right)\Im\left\{\upsilon\left(\tau_1\right)\right\}\right]
	\end{aligned},\label{eq34}
\end{equation}
where $\varepsilon\in\mathbb{R}$ fixes the size of the torus. And its frequencies are $\boldsymbol{\omega}_1^2=\left[2\pi/T_1,\alpha/T_1\right]^\mathrm{T}\in\mathbb{R}^{2\times1}$. $\Re\{\cdot\},\,\Im\{\cdot\}$ represent the notations of the real part and imagery part, respectively.

\section{Examples of application on nonlinear systems}\label{sec5}

In this section, all testing examples have been performed on a computer with 16 GB of RAM. The programming language is MATLAB and the version R2021b is used. And the flowchart of m-VCF is shown in Appendix \ref{appE}.

In the first example, the feasibility of the \textit{m}-VCF for the computation of the quasi-periodic solutions with multiple base frequencies is verified and validated. For the second example, the phase condition and Initialization technique in the \textit{m}-VCF are used to solve the quasi-periodic solutions with one unknown base frequency after a NS bifurcation of periodic solutions. The third example is used to verify the potential of \textit{m}-VCF for the computation of finite-element system.

\subsection{A single DOF Duffing-van der Pol oscillator}\label{sec5.1}

A \textit{d}-dimensional quasi-periodically forced Duffing-van der Pol oscillator is considered in this subsection:
\begin{equation}
	x-\mu\left(1-x^2\right)\dot{x}+\omega_0^2x+\alpha x^3=\sum_{i=1}^df_i\cos\left(\omega_it\right),\label{eq35}
\end{equation}
where $\mu$ represents positive real parameters, $\alpha$ is the coefficient controlling the nonlinearity, $\omega_0$ is modal frequency of corresponding linear system. $f_i$ and $\omega_i$ are the amplitude and frequency of the external excitation $f_i\cos\left(\omega_it\right)$. Based on this system, there are two practical examples computed to verify the presented methods:

In order to verify that the presented method can be used to solve the quasi-periodic solution with multiple dimensional base frequencies, Example 1 calculates the forced response of the system when \textit{d} is 1, 2 and 3 in external excitation $\sum_{i=1}^df_i\cos\left(\omega_it\right)$, respectively. Since the excitation frequencies are defined as incommensurable, the corresponding responses are periodic solutions $d=1$ and quasi-periodic solutions $d=2,\,3$. Specific parameters are as follows: $\mu=0.2,\,\alpha=0.5,\,\omega_0=2,\,f_1=2,\,f_2=1,\,f_3=0.5,\,\omega_1=\sqrt{5}\omega_2$ and $\omega_1=\sqrt{11}\omega_3$. Defining the first frequency $\omega_1$ as the continuation parameter, there are $d-1$ equations deficit in continuation. The so-called frequency conditions between  $\omega_1$ and $\omega_j,\,j=2,3$ are needed to add in Eqs. \eqref{eq9} and \eqref{eq10}.

\begin{figure*}[htbp]
	\centering
	\includegraphics[width=0.35\textwidth]{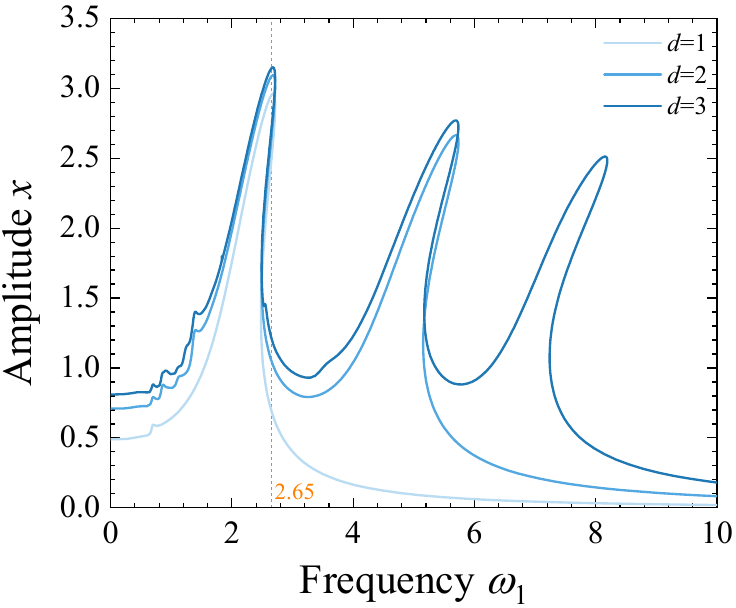}
	\caption{FRCs for the periodic (\textit{d}=1) and quasi-periodic (\textit{d}=2,3) orbits of the Duffing-van der Pol oscillator.}\label{fig1}
\end{figure*}

Fig. \ref{fig1} shows the FRCs of Duffing-van der Pol oscillator along the continuation parameter $\omega_1$. All FRCs in Fig. \ref{fig1} are calculated by the \textit{m}-VCF based on HB+…+HB, whose corresponding parameters are set as $u_i=[1,2,\cdots,5],\mathrm{~U}_i=11,\mathrm{~S}_i=2^5$, with $i=1,2,3$. These three curves are periodic (\textit{d}=1) and quasi-periodic (\textit{d}=2,3) orbits of system, respectively. Due to the increase in the number \textit{d} of the excitation frequencies, the corresponding FRC exhibits more complex dynamic behavior, where the FRC of \textit{d}-tori has \textit{d} main peaks. In addition, for peaks caused by the same frequency, the amplitude of the three quasi-periodic orbits is also nearly equal.

\begin{figure*}[htbp]
	\centering
	\includegraphics[width=0.85\textwidth]{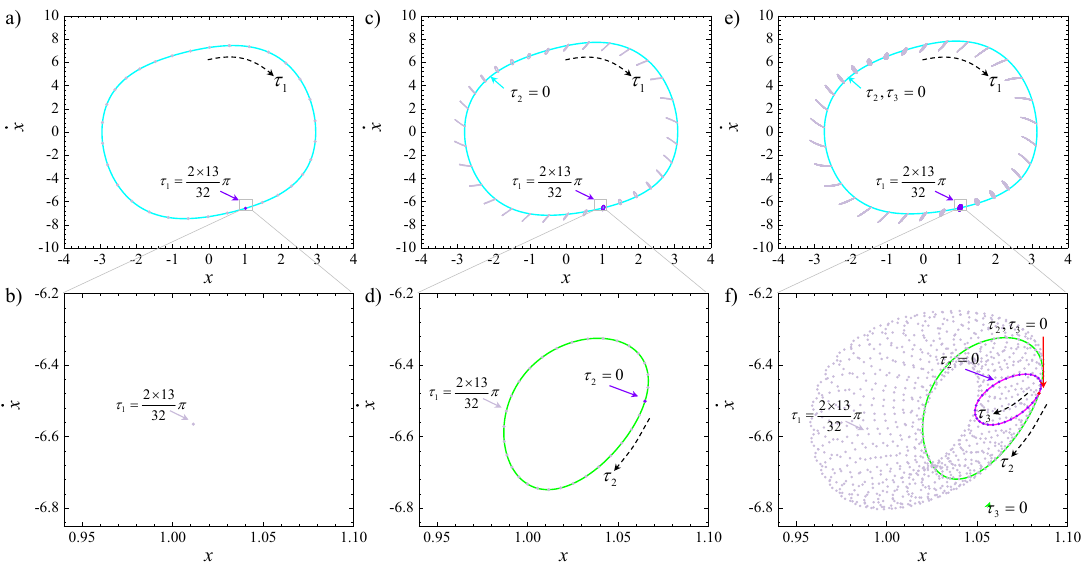}
	\caption{The quasi-periodic solution in the S domain for $\omega_1=2.65$, Phase diagrams: a) \textit{d}=1, c) \textit{d}=2, e) \textit{d}=3; Poincare sections: b) \textit{d}=1, d) \textit{d}=2, f) \textit{d}=3.}\label{fig2}
\end{figure*}

To reveal the relationship and difference of these three orbits, Fig. \ref{fig2} shows the phase diagrams and Poincare sections of these three orbits in the S domain for $\omega_{1}=2.65$. For phase diagram form Fig. \ref{fig2}a, \ref{fig2}c to \ref{fig2}e, the cycles marked by cyan, where $\tau_2,\tau_3$ are required to be equal to 0, have the same direction of $\tau_1$ and their shapes are almost similar. Choosing $\tau_{1}=0.8125\pi$, the Poincare sections are a point, a cycle and a torus, respectively. Similarly, when $\tau_3$ is required to be equal to 0, the cycles marked by green are also almost similar. Note that the cycles marked by cyan, green and magenta represent the periodicity of $\tau_i=\omega_it$, respectively. And the scale of the cycles with different colors decrease as the \textit{i} increases.

\begin{figure*}[htbp]
	\centering
	\includegraphics[width=0.85\textwidth]{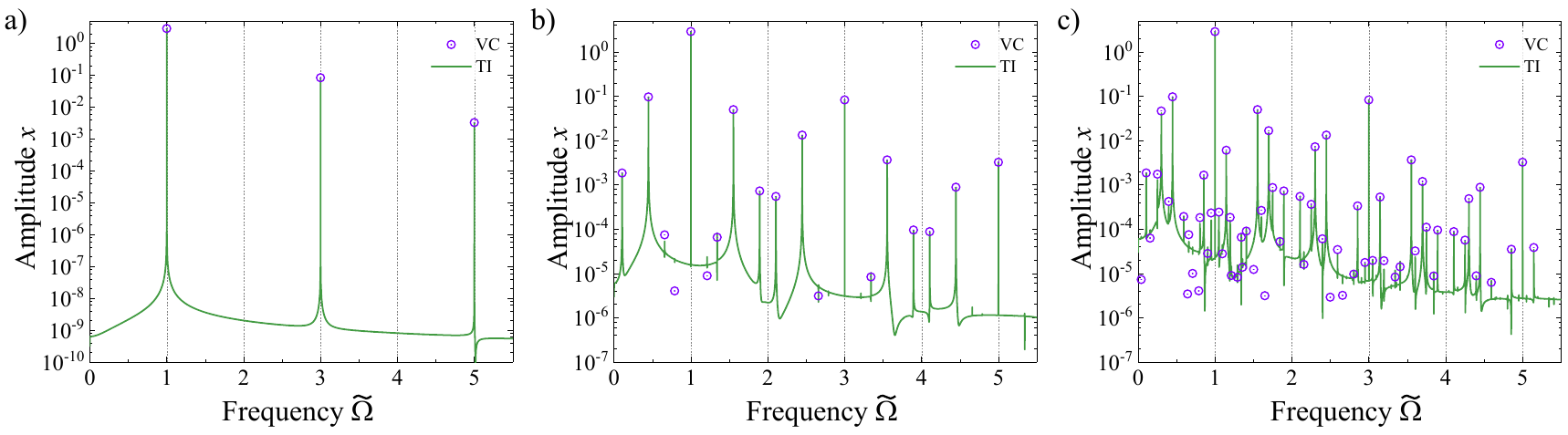}
	\caption{The comparisons of the quasi-periodic solutions computed by the \textit{presented method} and TI in the frequency domain: a) \textit{d}=1, b) \textit{d}=2, c) \textit{d}=3.}\label{fig3}
\end{figure*}

In addition, the comparisons of the quasi-periodic solutions computed by the presented method and TI in the frequency domain are also shown respectively in Fig. \ref{fig3}, where purple points represent the results by the presented method and green curves are the numerical solutions by TI and the notation $\tilde{\Omega}$ is set as $\omega/\omega_1$. The results of two methods are well consistent in Fig. \ref{fig3}. Moreover, as \textit{d} increases, the solutions of system become more complicated. Table \ref{tab1} concludes the results of these three FRCs by the presented method. The number $11^d+d$ of unknowns of nonlinear algebraic system \eqref{eq7} in the continuation increases exponentially as d increases. Correspondingly, more complex solutions take longer to compute.   

\begin{table*}[htbp]
	\centering
	\caption{The results of three FRCs of the Duffing-van der Pol oscillator.}\label{tab1}
	\small
	\begin{tabularx}{\textwidth}{XXXXXX}
		\toprule
		~ & Number of unknowns & Total time (s) & Number of points in FRC & Number of iterations in continuation & Time per point (s)  \\ 
		\midrule
		$d=1$ & 12 & 0.7 & 417 & 2451 & 0.0017  \\ 
		$d=2$ & 123 & 7.7 & 847 & 4760 & 0.0091  \\ 
		$d=3$ & 1334 & 2183.7 & 1295 & 7141 & 1.6863  \\
		\bottomrule
	\end{tabularx}
\end{table*}

\begin{figure*}[htbp]
	\centering
	\includegraphics[width=0.55\textwidth]{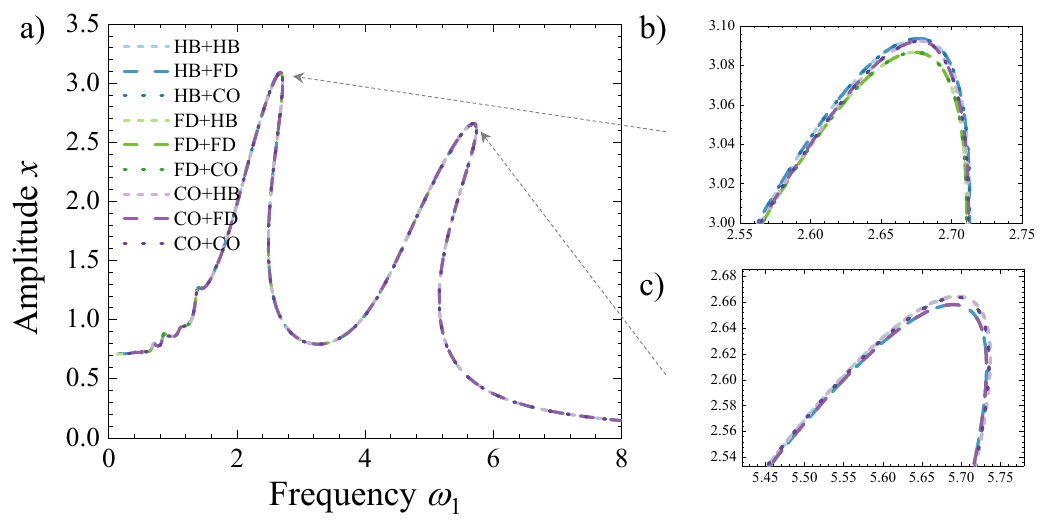}
	\caption{a): FRCs for the quasi-periodic (\textit{d}=2) orbits of the Duffing-van der Pol oscillator, b), c): Zoom of peaks.}\label{fig4}
\end{figure*}

Then the \textit{m}-VCF are used to calculate the quasi-periodic solutions with two base frequencies of Duffing-van der Pol oscillator under the external excitation $2\cos(\omega_1t)+\cos(\omega_2t)$ with $\omega_1=\sqrt{5}\omega_2$ by nine methods constructed by the \textit{m}-VCF, namely, HB + HB, HB + FD, HB + CO, FD + HB, FD + FD, FD + CO, CO + HB, CO + FD and CO + CO. Herein, the corresponding parameters of HB are set as $u_i=[1,2,\cdots,5],\,\mathrm{U}_i=11,\,\mathrm{S}_i=2^5$, the parameters of CO are $m=4,\,\mathrm{P}_i=8,\,\mathrm{U}_i=2^5,\,\mathrm{S}_i=2^5$, and the parameters of FD are $K=\left[-3,-2,-1,0,1\right],\,\mathrm{U}_i=2^5,\,\mathrm{S}_i=2^5$. Fig. \ref{fig4} shows the FRCs for the quasi-periodic (\textit{d}=2) orbits of Duffing-van der Pol oscillator along the continuation parameter $\omega_1$. The results show that all nine constructed methods can compute the solutions of Duffing-van der Pol oscillator with good consistency. The error in the peak is mainly due to the fact that the $\mathrm{U}_{i}$ of the three methods are not very large, and it does not play a good approximation at the peaks. However, the results are sufficient to demonstrate the feasibility of the presented method. In addition, Table \ref{tab2} concludes the results of the FRCs for the quasi-periodic (\textit{d}=2) orbits by the nine constructed methods of the \textit{m}-VCF. Note that since HB+HB has the fewest unknowns, the calculation time is also the smallest. 

\begin{table*}[htbp]
	\centering
	\caption{The results of the FRCs for the quasi-periodic (\textit{d}=2) orbits of the Duffing-van der Pol oscillator by the nine constructed methods of the \textit{m}-VCF.}\label{tab2}
	\small
	\begin{tabularx}{\textwidth}{XXXXXX}
		\toprule
		~ & Number of unknowns & Total time (s) & Number of points in FRC & Number of iterations in continuation & Time per point (s)  \\ 
		\midrule
		HB+HB & 123 & 7.7 & 847 & 4760 & 0.0091     \\
		HB+FD & 354 & 35.3 & 809 & 5117 & 0.0436    \\
		HB+CO & 354 & 39.3 & 806 & 5168 & 0.0488    \\
		FD+HB & 354 & 59.0 & 839 & 5220 & 0.0703    \\
		FD+FD & 1026 & 863.4 & 963 & 6690 & 0.8965  \\
		FD+CO & 1026 & 910.8 & 964 & 6694 & 0.9448  \\
		CO+HB & 354 & 61.2 & 847 & 5382 & 0.0722    \\
		CO+FD & 1026 & 1044.3 & 980 & 6753 & 1.0656 \\
		CO+CO & 1026 & 957.1 & 972 & 6781 & 0.9847 \\
		\bottomrule
	\end{tabularx}
\end{table*}

\subsection{A cantilevered pipe system}\label{sec5.2}

\begin{figure*}[htbp]
	\centering
	\includegraphics[width=0.4\textwidth]{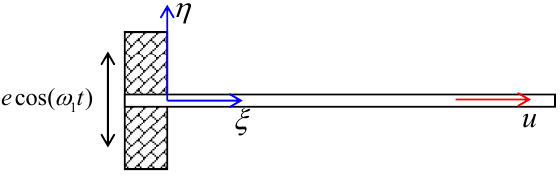}
	\caption{A cantilevered pipe conveying fluid subjected to a periodic excitation.}\label{fig5}
\end{figure*}

The second system is a conveying fluid conveying fluid subject to a periodic excitation, as shown in Fig. \ref{fig5} mentioned in \cite{JW40}, where $\eta$ and $\xi$ are dimensionless transverse displacement and arc-length, $u$ is flow velocity.

Using \textit{n}-mode Galerkin discretization $\eta(\xi,t)=\sum_{j=1}^{n}\phi_{j}(\xi)q_{j}(t)$, the cantilevered pipe system can be governed by:
\begin{equation}
	\begin{aligned}
		&m_{ij}\ddot{q}_{j}+c_{ij}\dot{q}_{j}+k_{\bar{i}j}q_{j}+\alpha_{ijkl}q_{j}q_{k}q_{l}
		\\&+\beta_{ijkl}q_jq_k\dot{q}_l+\gamma_{ijkl}q_j\dot{q}_k\dot{q}_l=g_ie\omega_1^2\cos\omega_1t
	\end{aligned},\label{eq36}
\end{equation}
where the coefficients $m_{ij},\,c_{ij},\,k_{ij}q_{j},\,\alpha_{ijkl},\,\beta_{ijkl},\,\gamma_{ijkl},\,g_{i}$ are computed from integrals of the eigenfunctions in \cite{JW40}. In this work, the mode of \textit{n} is chosen as 4.

\begin{figure*}[htbp]
	\centering
	\includegraphics[width=0.75\textwidth]{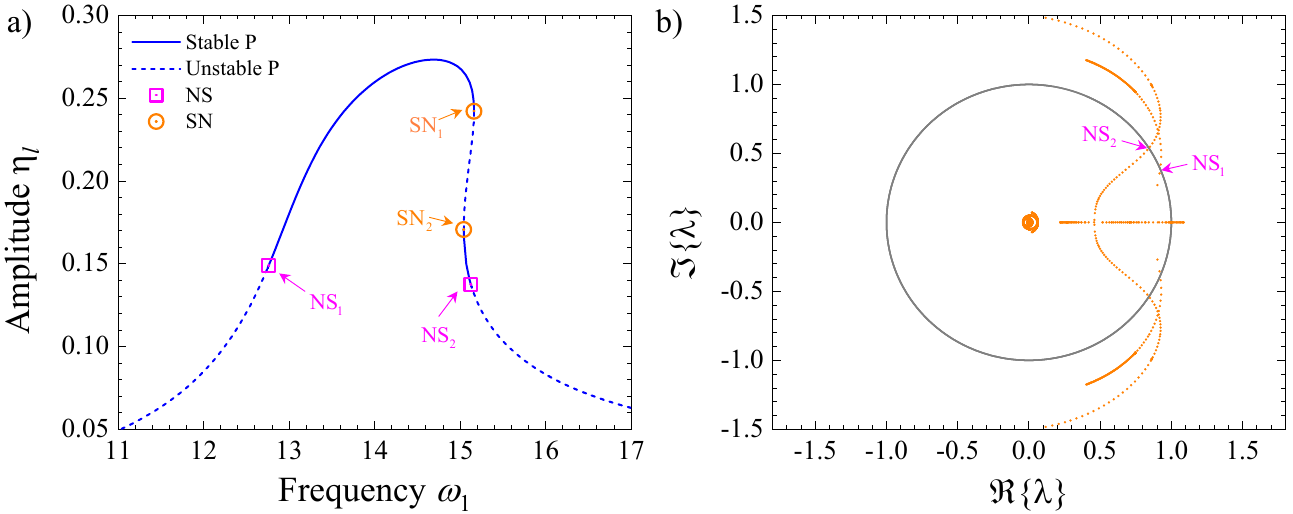}
	\caption{a) FRC for the periodic (\textit{d}=1) orbits of the cantilevered pipe system with $e=0.012$, b) The eigenvalues of monodromy matrix.}\label{fig6}
\end{figure*}

In Fig. \ref{fig6}a, the FRC for the periodic orbits of the cantilevered pipe system with $e=0.012$ is calculated by the \textit{m}-VCF constructed by the FD, where the parameters are set as $K=\left[-3,-2,-1,0,1,2\right],\,\mathrm{U}_{1}=55,\,\mathrm{S}_{1}=55$. Here, the blue solid curves represent stable solutions and the blue dotted curves stand for unstable solutions, the squares represent the Neimark-Sacker bifurcations (NS), and circles are saddle–node bifurcations (SN). The amplitude $\eta_{l}$ is the maximum of the dimensionless transverse displacement at the free end. The stabilities and bifurcations of the periodic solutions are assessed by the Floquet method in Appendix \ref{appD}, and the eigenvalues of monodromy matrix are shown in Fig. \ref{fig6}b. At a NS bifurcation, the quasi-periodic solution (\textit{d}=2) occurs when the periodic solutions are unstable.

\begin{figure*}[htbp]
	\centering
	\includegraphics[width=0.75\textwidth]{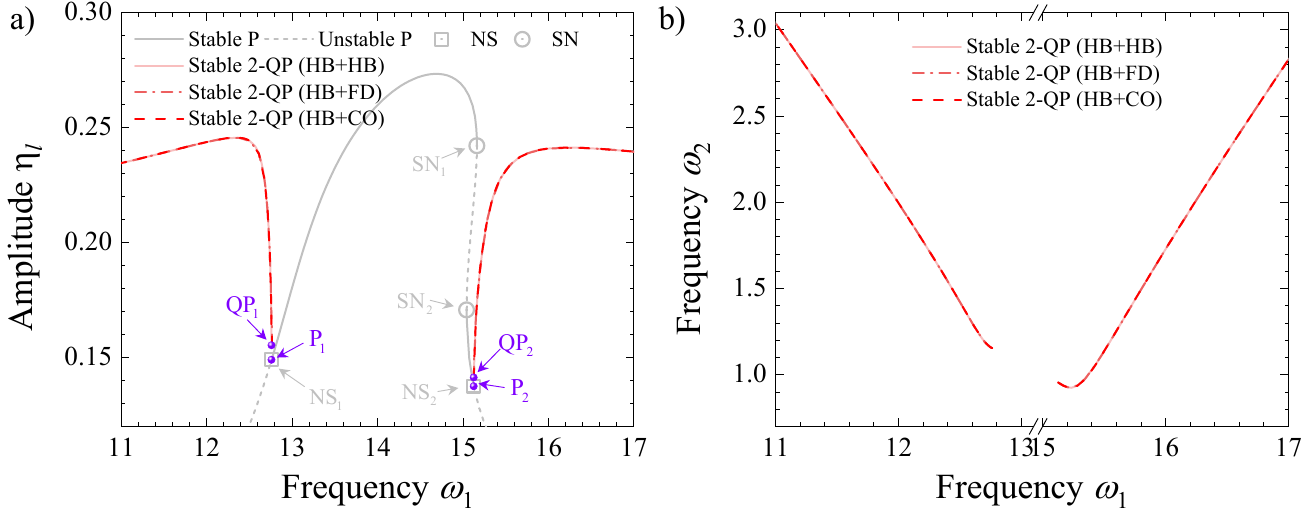}
	\caption{a) FRC for the quasi-periodic periodic (d=2) orbits of the cantilevered pipe system with $e=0.012$, b) The relationship between the new frequency $\omega_2$ and frequency $\omega_1$ of excitation.}\label{fig7}
\end{figure*}

Then, to demonstrate the feasibility of phase condition $PC_{2}$ in subsection \ref{sec4.1}, the function of $\boldsymbol{\Phi}_2(\tau_2)$ is approximated by HB, FD and CO, respectively. Here, the corresponding parameters of HB are set as $u_2=[1,2,\cdots,7],\mathrm{~U}_2=15,\mathrm{~S}_2=2^6$, the parameters of CO are $m=4,\mathrm{~P}_2=16,\mathrm{~U}_2=2^6,\mathrm{~S}_2=2^6$ and the parameters of FD are $K=\left[-3,-2,-1,0,1,2\right],\mathrm{~U}_2=2^6,\mathrm{~S}_2=2^6$. And $\boldsymbol{\Phi}_1(\tau_1)$ is approximated by HB with $k^1=[1,2,\cdots,7],\mathrm{~U}_1=15,\mathrm{~S}_1=2^6$. The FRCs and the relationship between two frequencies for the quasi-periodic (\textit{d}=2) orbits of the cantilevered pipe system are shown in Fig. \ref{fig7}, where the methods constructed by HB+HB, HB+FD and HB+CO are used. The results marked by three kind of red curves are well consistent in Fig. \ref{fig7}. And the stabilities of the quasi-periodic (\textit{d}=2) solutions are assessed by the Lyapunov exponents in Appendix \ref{appD}. The evolution of the first four Lyapunov exponents of the quasi-periodic solutions calculated by HB+HB are shown in Fig. \ref{fig8}, where the parameters in Appendix \ref{appD} are set as $j=1,\,Y_2=2^6,\,N_\mathrm{L}=10^4,\,N_\mathrm{M}=2^8$. It can be seen that the exponents of the quasi-periodic solutions are all not greater than 0. That is to say, they are all stable. The results of HB+FD and HB+CO are similar and will not be shown in detail here. To analyze the different motions of the pipe deformation along the pipe central axis, Fig. \ref{New fig. 2} shows that the results of unstable periodic solutions and quasi-periodic solutions in 40 periods $T_1=2\pi/\omega_1$ at $\omega_1=11,\,17$, where red and blue curves represents the trajectories of transverse displacement at $\xi=l,\,l/2$. Quasi-periodic solutions exhibit complex dynamical behaviors, and their amplitudes are much larger than those of the periodic solutions.

\begin{figure*}[htbp]
	\centering
	\includegraphics[width=0.75\textwidth]{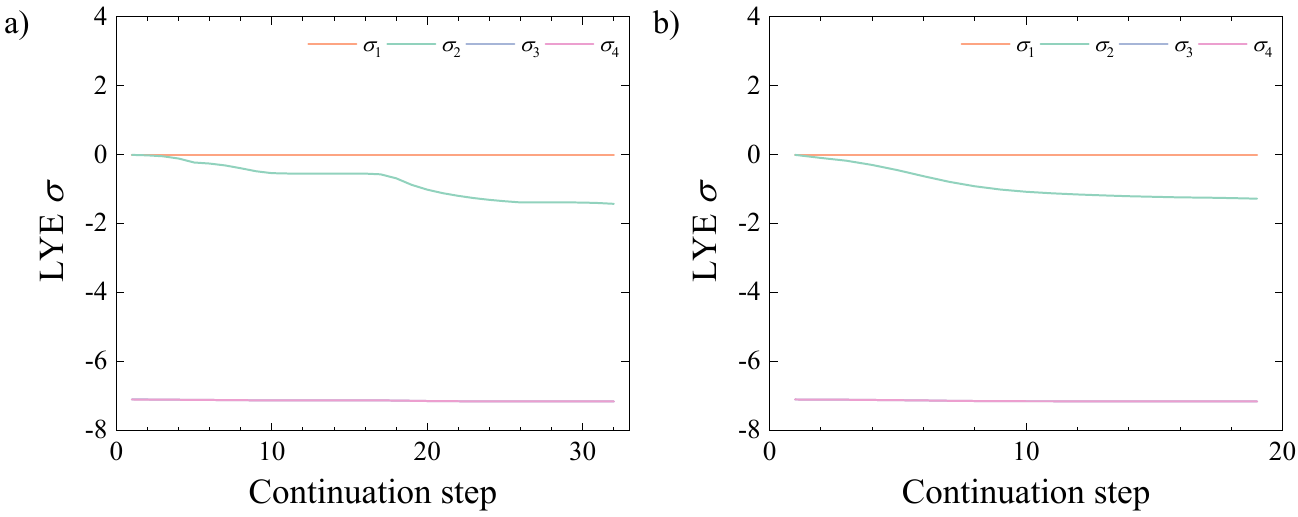}
	\caption{The evolution of the first four Lyapunov exponents of the quasi-periodic solutions of the cantilevered pipe system calculated by HB+HB: a) The FRC of QP on the left, b) The FRC of QP on the right.}\label{fig8}
\end{figure*}

\begin{figure*}[htbp]
	\centering
	\includegraphics[width=0.75\textwidth]{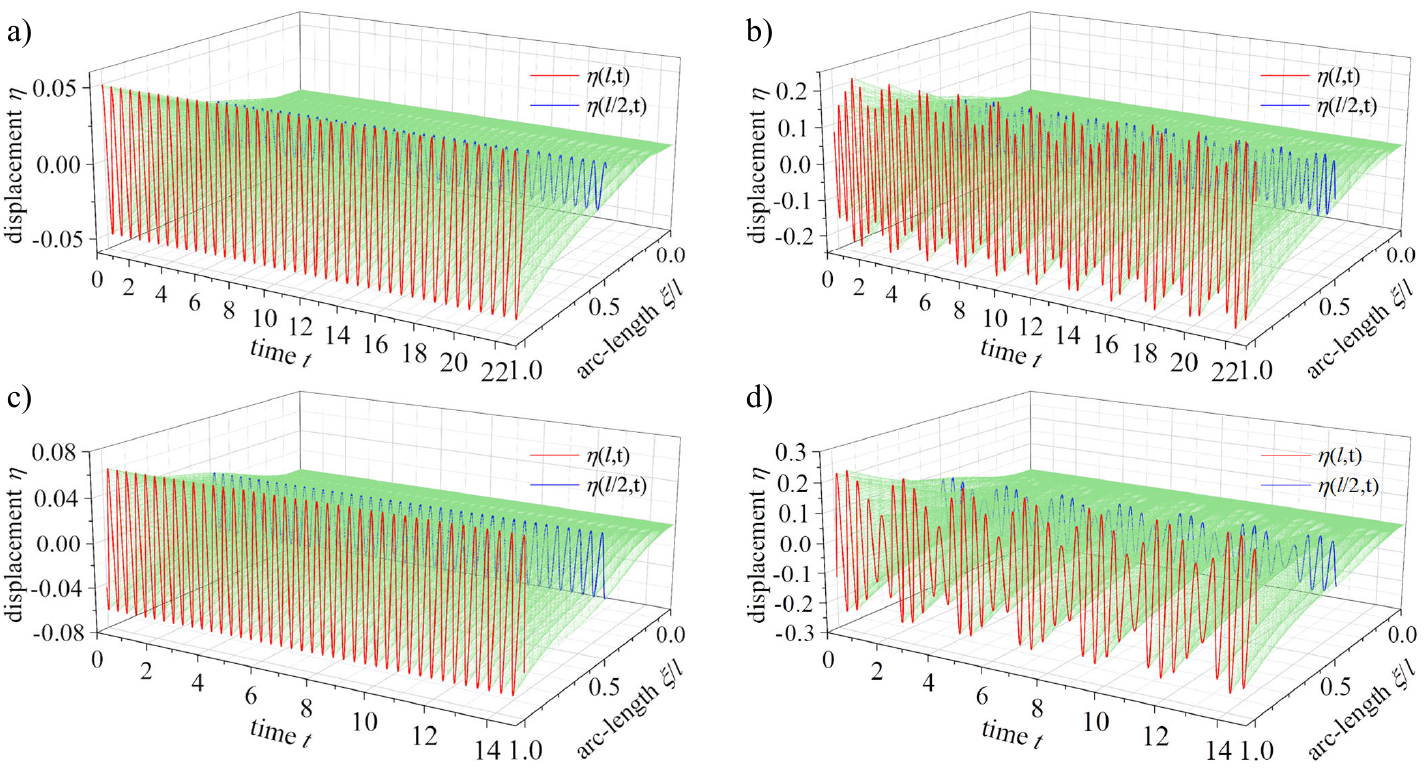}
	\caption{The transverse displacement of pipe along the pipe central axis: a) Unstable periodic solution at $\omega_1=11$, b) Stable quasi-periodic solution at $\omega_1=11$, c) Unstable periodic solution at $\omega_1=17$, d) Stable quasi-periodic solution at $\omega_1=17$.}\label{New fig. 2}
\end{figure*}

In addition, Table \ref{tab3} concludes the results of the FRCs for the quasi-periodic (\textit{d}=2) orbits. It should be noted that the time in Table \ref{tab3} is the total time to solve the quasi-periodic solutions and evaluate their stabilities. It is not difficult to find that HB+HB has the fewest unknowns and the smallest calculation cost. And the time used to analyze the stability is almost equivalent, because the steps used to assess the stability of the solutions calculated by the three methods, HB+HB, HB+FD and HB+CO, are the same.

\begin{table*}[htbp]
	\centering
	\caption{The results of FRCs for the quasi-periodic (\textit{d}=2) orbits on the left of the cantilevered pipe system.}\label{tab3}
	\small
	\begin{tabularx}{\textwidth}{XXXXXXXX}
		\toprule
		~ & Number of unknowns & Total time of computation (s) & Total time of evaluation (s) & Number of points in FRC & Number of iterations in continuation & Time per point of computation (s) & Time per point of evaluation (s) \\ 
		\midrule
		HB+HB & 902 & 24.6 & 24.8 & 32 & 183 & 0.7688 & 0.7750  \\ 
		HB+FD & 3842 & 1140.7 & 42.1 & 55 & 415 & 20.7400 & 0.7655  \\ 
		HB+CO & 3842 & 1289.3 & 49.2 & 63 & 401 & 20.4651 & 0.7810  \\
		\bottomrule
	\end{tabularx}
\end{table*}

\begin{figure*}[htbp]
	\centering
	\includegraphics[width=0.75\textwidth]{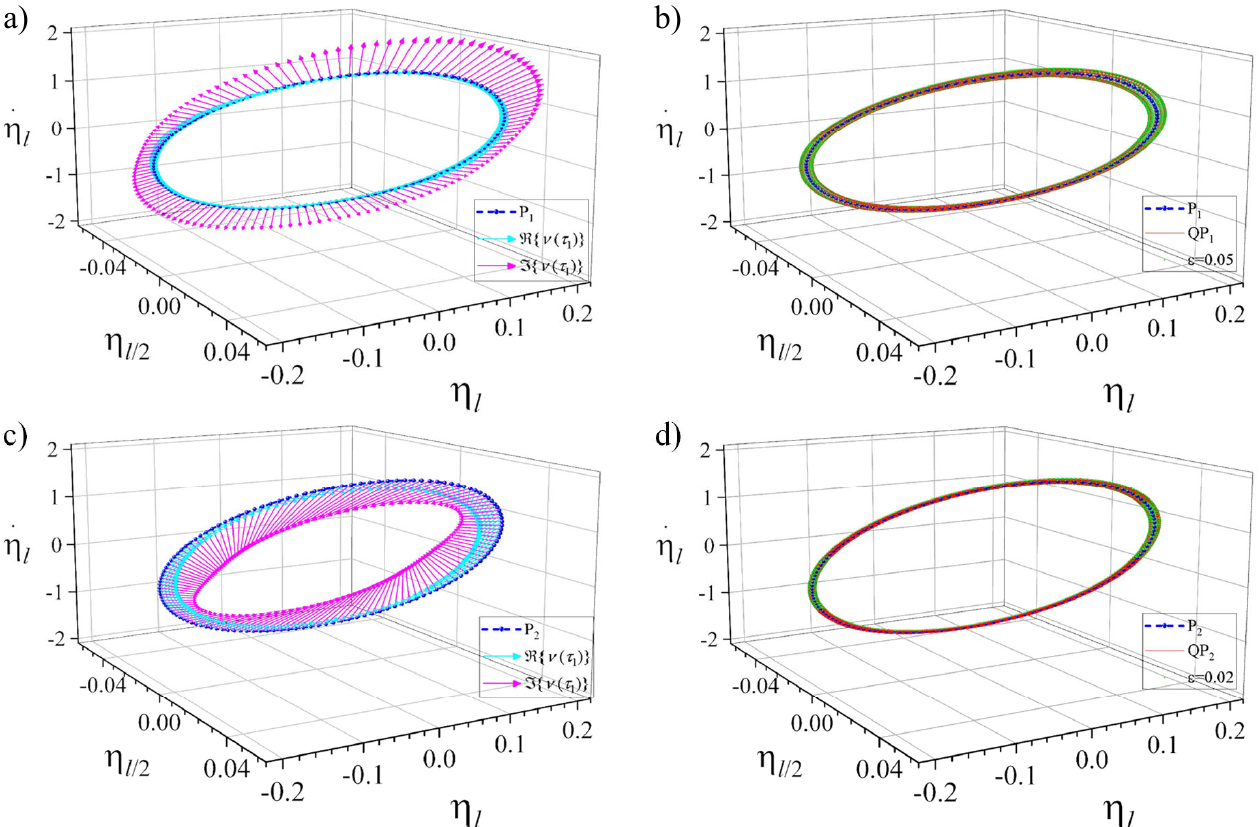}
	\caption{Initialization of quasi-periodic solution of $\mathrm{QP}_1$ and $\mathrm{QP}_2$ based on $\mathrm{P}_1$ and $\mathrm{P}_2$ of the cantilevered pipe system: a) and c), The eigenvectors $\upsilon(\tau_{1})$ of the monodromy matrix; b) and d), The comparison between the solution of initialization and the reference solution calculated by HB+HB.}\label{fig9}
\end{figure*}

Note that the quasi-periodic solutions of $\mathrm{QP}_1$ and $\mathrm{QP}_2$ are firstly predicted by the periodic solution of $\mathrm{P}_1$ and $\mathrm{P}_2$ by means of the initialization of quasi-periodic solution with two base frequencies in subsection \ref{sec4.2}. Then the solutions of $\mathrm{QP}_1$ and $\mathrm{QP}_2$ are determined by the method constructed by HB+HB. Fig. \ref{fig9} shows the initialization of quasi-periodic solution of $\mathrm{QP}_1$ and $\mathrm{QP}_2$ based on $\mathrm{P}_1$ and $\mathrm{P}_2$. Here, the eigenvectors $\upsilon(\tau_{1})$ of the monodromy matrix $\boldsymbol{M}(\tau_1+2\pi,\tau_1)$ of $\mathrm{P}_1$ and $\mathrm{P}_2$ are shown in Fig. \ref{fig9}a and Fig. \ref{fig9}c. Defining $\varepsilon=0.05$ and $\varepsilon=0.02$, the $\mathrm{QP}_1$ and $\mathrm{QP}_2$ marked by green points are initialized in Fig. \ref{fig9}b and Fig. \ref{fig9}d, respectively. The red curves are the reference solution calculated by HB+HB. The results in Fig. \ref{fig9} shows that the initialization in subsection \ref{sec4.2} is a good predictor of quasi-periodic solutions after NS bifurcation. 

\subsection{A Bernoulli beam with nonlinear support spring}\label{sec5.3}

The third system is a clamped-free Bernoulli beam with a linear and nonlinear springs support at its free end \cite{JW52}, as shown in Fig. \ref{new fig. 1}, whose force-displacement relation of the nonlinear force is given by:
\begin{figure*}[htbp]
	\centering
	\includegraphics[width=0.7\textwidth]{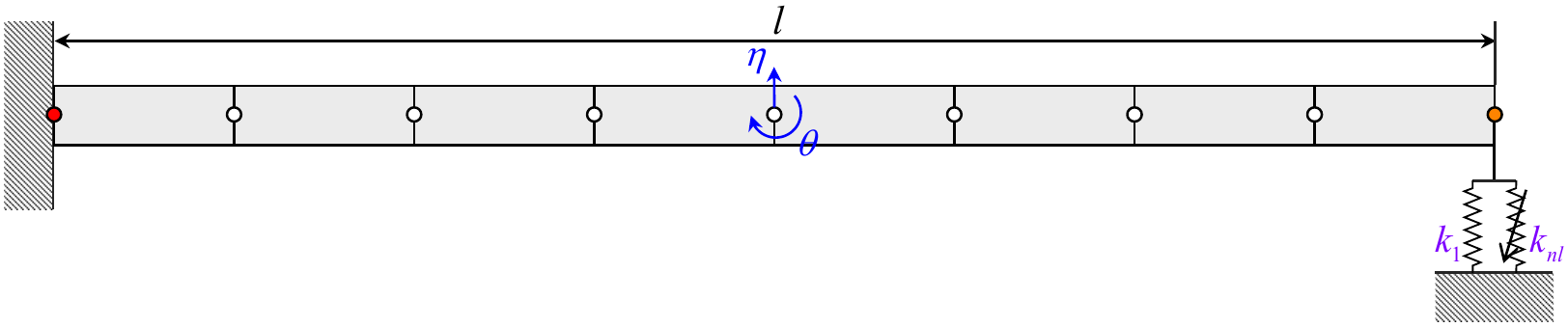}
	\caption{A Finite-element clamped-free Bernoulli beam with a linear and nonlinear cubic support springs at its free end.}\label{new fig. 1}
\end{figure*}

\begin{equation}
	F=k_1\eta+k_{nl}\eta^3,\label{eq37}
\end{equation}
where $\eta$ is the transverse displacement, $F$ is the spring force and $k_1=27\mathrm{~N/mm}$, $k_{nl}=60\mathrm{~N/mm^3}$ denote linear and nonlinear stiffness, respectively.

The length, width and height of the beam are set as $l=2700\mathrm{mm},\,h=b=10\mathrm{mm}$, density is $\rho=1780{\times}10^{-9}\mathrm{~kg/mm^3}$ and Young’s modulus is $E=45\times10^6\mathrm{kPa}$. Using classic finite-element method, the transverse displacement $\eta$ and the rotation angle $\theta$ are introduced at each node. Define $N_{e}$ is the number of elements used in the finite-element method, there are $2N_{e}$ degrees of freedom since the transverse displacement and the rotation angle of the clamped end are always 0. The cantilever Bernoulli beam can be modeled as:
\begin{equation}
	\boldsymbol{M}\ddot{\boldsymbol{x}}+\boldsymbol{D}\dot{\boldsymbol{x}}+\boldsymbol{K}\boldsymbol{x}+\boldsymbol{F}_{nl}(\boldsymbol{x})=e\boldsymbol{f}\cos\omega_1t,\label{eq38}
\end{equation}
where $\boldsymbol{M},\,\boldsymbol{D},\,\boldsymbol{K}$ are the mass, the damping and the stiffness matrices, $\boldsymbol{F}_{nl}(x)$is nonlinear forces introduced by the nonlinear stiffness at the free end. And the damping matrix is set as $\boldsymbol{D}=\alpha\boldsymbol{M}+\boldsymbol{K}_b$ with $\alpha=1.25\times10^{-4}\mathrm{s}$ and $\beta=2.5\times10^{-5}\mathrm{s}^{-1}$. The stiffness matrix of beam elements $\boldsymbol{K}_b$ is obtained from $\boldsymbol{K}$ by removing the contribution of the linear stiffness $k_1$ of the support spring. And the vector of $\boldsymbol{f}$ is set as $f=\omega_{l_{1}}^{2}\boldsymbol{M\phi}_{\mathrm{l}}$ with $\omega_{l_{1}}=15.60\mathrm{rad/s}$ being the first natural frequency of the undamped linear system and $\boldsymbol{\phi}$ being the linear normal mode corresponding to $\omega_{l_{1}}$.

\begin{figure*}[htbp]
	\centering
	\includegraphics[width=0.75\textwidth]{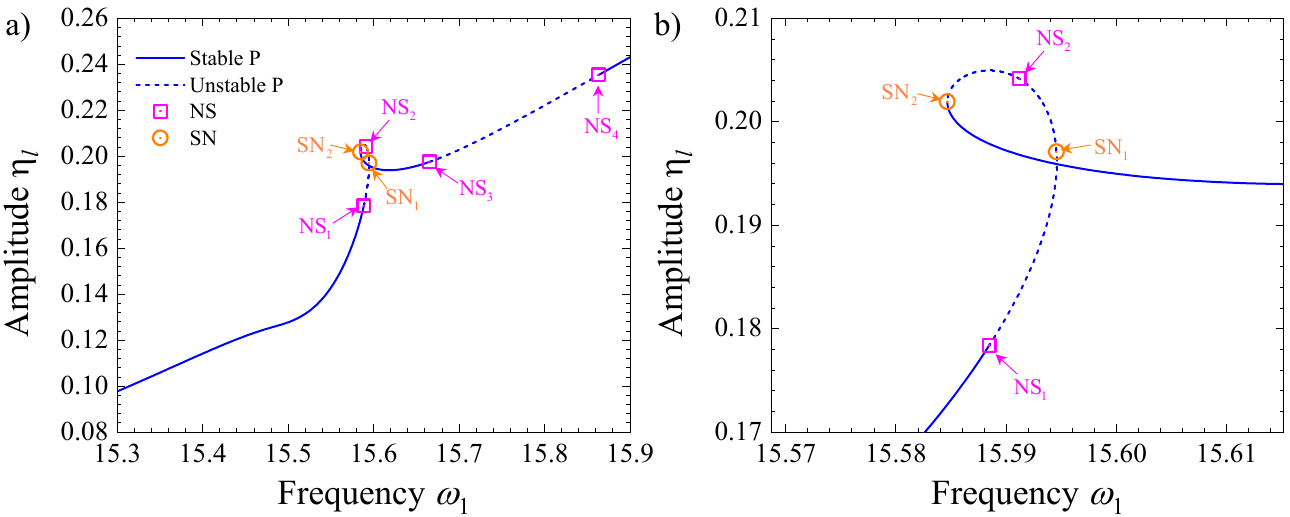}
	\caption{FRCs for the periodic (\textit{d}=1) orbits of the Bernoulli beam, b): Zoom of peaks.}\label{fig10}
\end{figure*}

With $e=0.002,\,N_e=8$, the FRC for the periodic orbits of the Bernoulli beam is shown in Fig. \ref{fig10}, where the blue solid curves represent stable solutions and the blue dotted curves stand for unstable solutions, the squares represent the NS bifurcations and circles are SN bifurcations. Here, the periodic solutions are calculated by the HB with $u_1=\left[1,2,\cdots,7\right],\,\mathrm{U}_1=15,\,\mathrm{S}_1=2^6$. The quasi-periodic solutions occur at the NS bifurcation points based on Floquet theory. 

\begin{figure*}[htbp]
	\centering
	\includegraphics[width=0.75\textwidth]{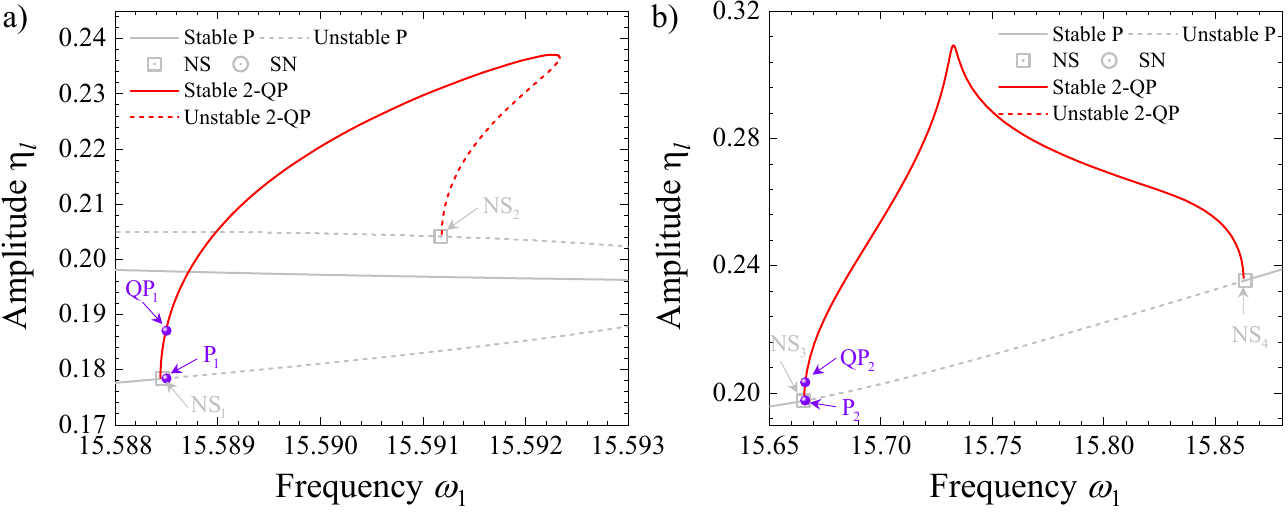}
	\caption{FRC for the quasi-periodic periodic (d=2) orbits of the Bernoulli beam, a): between the $\mathrm{NS}_1$ and $\mathrm{NS}_2$ bifurcations; b): between the $\mathrm{NS}_3$ and $\mathrm{NS}_4$ bifurcations.}\label{fig11}
\end{figure*}

\begin{figure*}[htbp]
	\centering
	\includegraphics[width=0.75\textwidth]{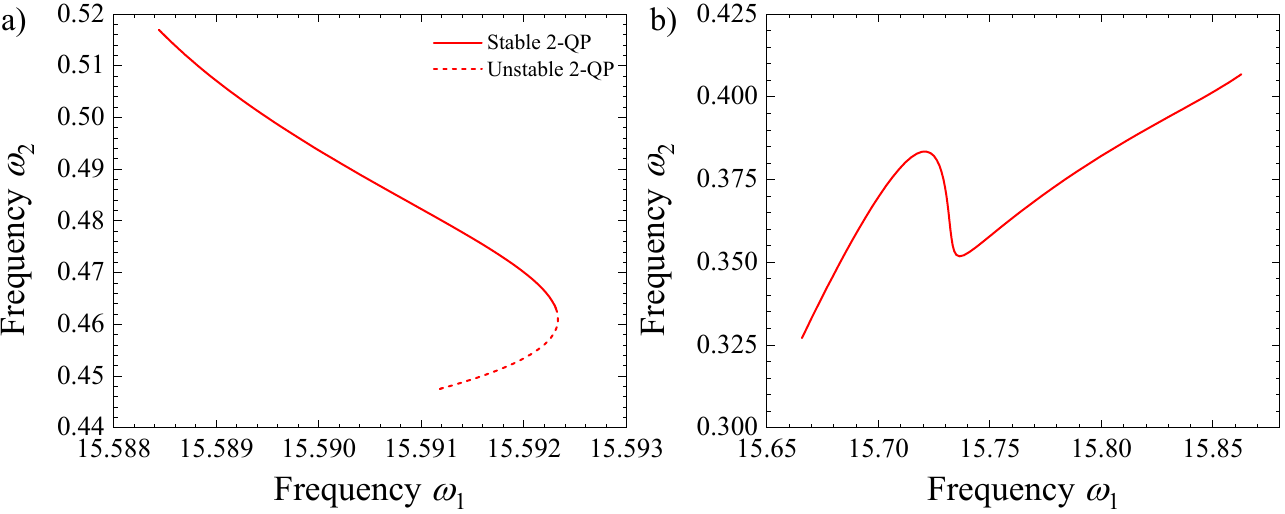}
	\caption{The relationship between the new frequency $\omega_2$ and frequency $\omega_1$ of excitation of the Bernoulli beam, a): between the $\mathrm{NS}_1$ and $\mathrm{NS}_2$ bifurcations; b): between the $\mathrm{NS}_3$ and $\mathrm{NS}_4$ bifurcations.}\label{fig12}
\end{figure*}

Based on the results of the cantilevered pipe system, HB+HB are used to save computing cost. The FRC for the quasi-periodic (\textit{d}=2) orbits of the Bernoulli beam is shown in Fig. \ref{fig11}, where the corresponding parameters of HB are set as $u_i=[1,2,\cdots,11],\,\mathrm{U}_i=23,\,\mathrm{S}_i=2^6$. Here, the red solid curves represent stable solutions and the red dotted curves stand for unstable solutions. The stabilities of quasi-periodic solutions in Fig. \ref{fig11}a change at $\omega_{1}=15.5923$, it can be considered as a SN bifurcation point. And the quasi-periodic solutions are all stable in Fig. \ref{fig11}b. The relationship between the new frequency $\omega_2$ and frequency $\omega_1$ of excitation of the quasi-periodic solutions are also shown in Fig. \ref{fig12}. The results showed that the new frequency $\omega_2$ is in the range of 0.32 to 0.52 and is smaller than the frequency $\omega_1$ of excitation. 

\begin{figure*}[htbp]
	\centering
	\includegraphics[width=0.75\textwidth]{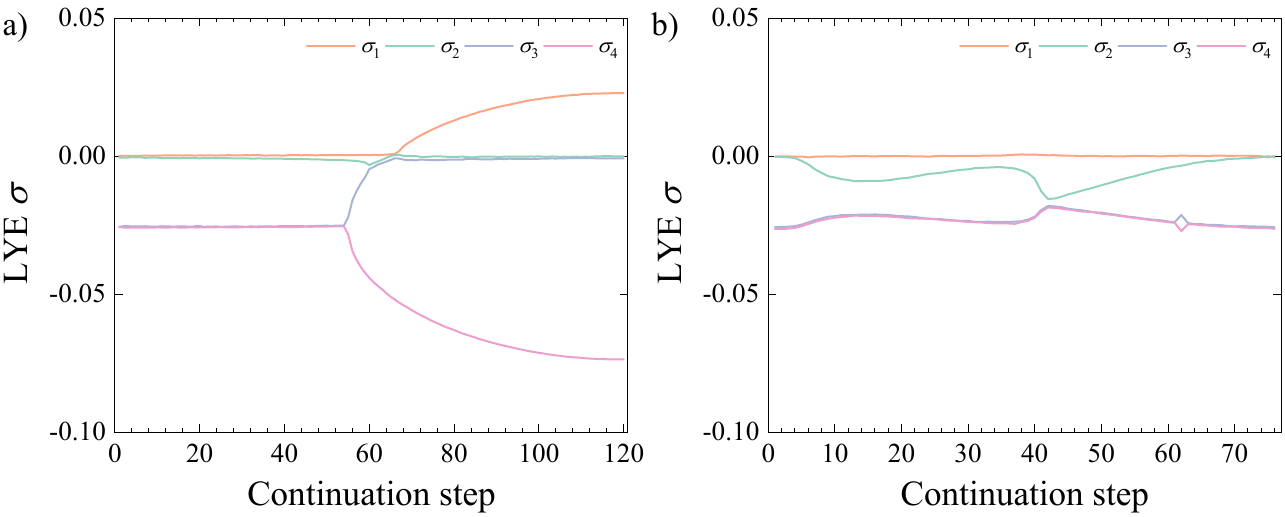}
	\caption{The evolution of the first four Lyapunov exponents of the quasi-periodic solutions of the Bernoulli beam, a): between the $\mathrm{NS}_1$ and $\mathrm{NS}_2$ bifurcations; b): between the $\mathrm{NS}_3$ and $\mathrm{NS}_4$ bifurcations.}\label{fig13}
\end{figure*}

\begin{figure*}[htbp]
	\centering
	\includegraphics[width=0.75\textwidth]{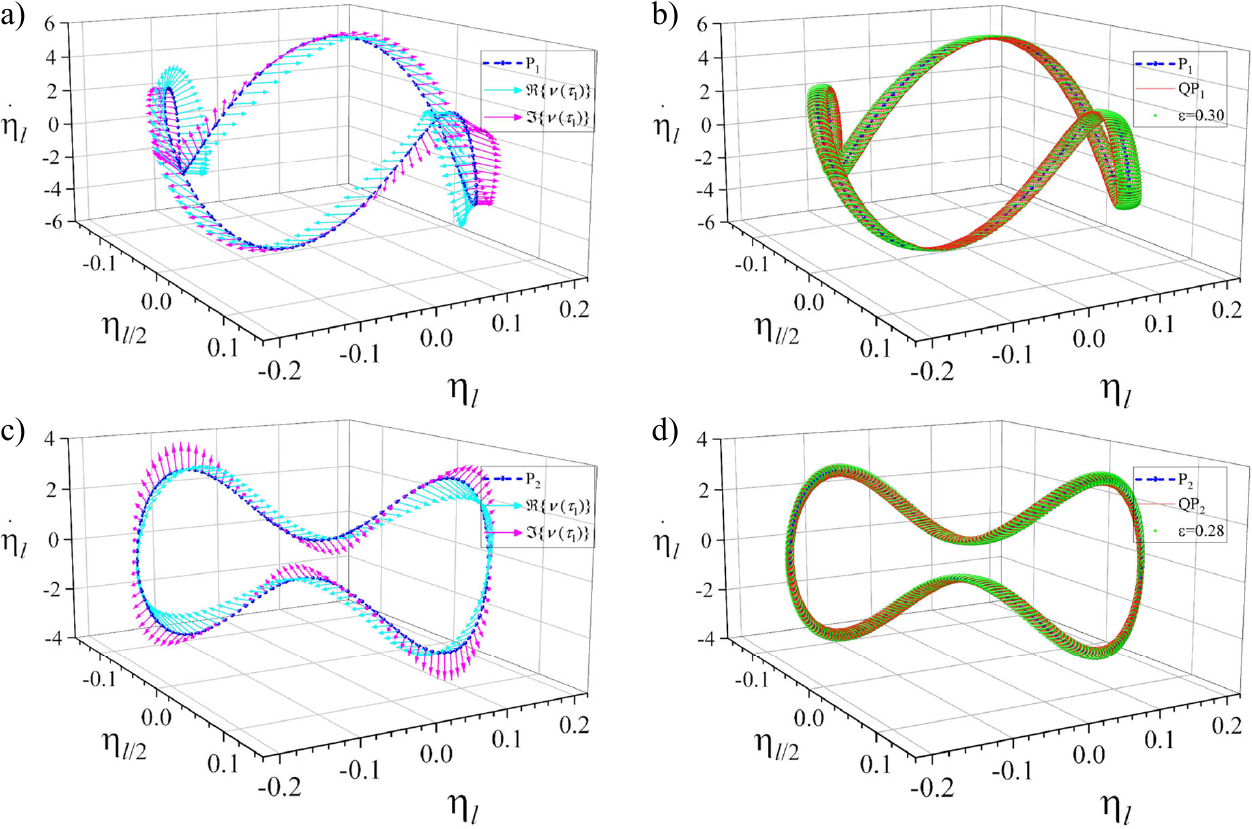}
	\caption{Initialization of quasi-periodic solutions of $\mathrm{QP}_1$ and $\mathrm{QP}_2$ based on $\mathrm{P}_1$ and $\mathrm{P}_2$ of the Bernoulli beam: a) and c), The eigenvectors $\upsilon(\tau_{1})$ of the monodromy matrix; b) and d), The comparison between the solution of initialization and the reference solution calculated by HB+HB.}\label{fig14}
\end{figure*}

Fig. \ref{fig13} shows the evolutions of the first four Lyapunov exponents of the quasi-periodic solutions, where the where the parameters in Appendix \ref{appD} are set as $j=1,\,Y_{2}=2^{6},\,N_{\mathrm{L}}=3\times10^{4},\,N_{\mathrm{M}}=2^{11}$. Table \ref{tab4} concludes the results of the FRCs for the quasi-periodic (\textit{d}=2) orbits, where Orbit 1 and Orbit 2 represent the quasi-periodic orbits of $\mathrm{NS}_1$-$\mathrm{NS}_2$ and $\mathrm{NS}_3$-$\mathrm{NS}_4$, respectively. It should be noted that the time in Table \ref{tab4} is the total time to solve the quasi-periodic solutions and evaluate their stabilities. Note that compared to the cantilevered pipe system, the time to evaluate stability increases due to the increase of $N_{\mathrm{L}}$ and $N_{\mathrm{M}}$. Similar with the cantilevered pipe system in subsection \ref{sec5.2}, Fig. \ref{fig14} shows the initialization of quasi-periodic solutions of $\mathrm{QP}_1$ and $\mathrm{QP}_2$ based on $\mathrm{P}_1$ and $\mathrm{P}_2$, where $\varepsilon=0.30$ and $\varepsilon=0.28$, respectively. The results show again that the initialization in subsection \ref{sec4.2} is a good predictor of quasi-periodic solutions after NS bifurcation. 

\begin{table*}[htbp]
	\centering
	\caption{The results of FRCs for the quasi-periodic (d=2) orbits of the Bernoulli beam.}\label{tab4}
	\small
	\begin{tabularx}{\textwidth}{XXXXXXXX}
		\toprule
		~ & Number of unknowns & Total time of computation (s) & Total time of evaluation (s) & Number of points in FRC & Number of iterations in continuation & Time per point of computation (s) & Time per point of evaluation (s) \\ 
		\midrule
		Orbit 1 & 8466 & 5024.9 & 1891.3 & 120 & 366 & 41.8742 & 15.7608 \\
		Orbit 2 & 8466 & 3592.3 & 1190.4 & 76 & 315 & 47.2671 & 15.6632  \\
		\bottomrule
	\end{tabularx}
\end{table*}

\section{Conclusions}\label{sec6}

In this work, a multi-steps variable-coefficient formulation (\textit{m}-VCF) is proposed based on a unified framework governed by Eq. \eqref{eq5}, where the matrices $\boldsymbol{\Upsilon}_i{}^0(k_i),\boldsymbol{\Upsilon}_i{}^1(k_i),\boldsymbol{\Upsilon}_i{}^2(k_i):\mathbb{U}_i\to\mathbb{R}^{\mathrm{U}_i\times\mathrm{U}_i}$ are constant because $k_i$ is a user-defined parameter based on either HB or CO or FD. Based on this unified framework, one can iteratively obtain the nonlinear algebraic equations of Eq. \eqref{eq7} for the quasi-periodic solutions with multiple base frequencies. Specifically, if the \textit{m}-VCF is implemented by HB in each step, it is actually the extension of VCHB for multiple base frequencies, called as \textit{m}-VCHB. Also, there are \textit{m}-VCFD and \textit{m}-VCCO. Moreover, the present \textit{m}-VCF rectifies the procedure for handling nonlinear force terms in the previous VCHB for the quasi-periodic solutions with two base frequencies because it is incorrectly assumed that $\mathbf{F}^1$ was only related to $\mathbf{Z}^1$. Furthermore, inspired by the AFT, an efficient method of alternating U and S domain is proposed in this work to evaluate the nonlinear force terms under \textit{m}-VCF. If the nonlinear force terms $\mathbf{F}^0$ in Eq. \eqref{eq1} corresponds only to $\mathbf{Z}^0$, \textcolor{red}{ the number of the multiplications for solving quasi-periodic solutions with d-base frequencies is $N_{m\mathrm{-VCHB}}=n^{2}\sum_{i=1}^{d}(\mathrm{S}_{i}\mathrm{U}_{i}+\mathrm{U}_{i}\mathrm{S}_{i}\mathrm{U}_{i})\mathrm{U}_{1}^{i-1}\mathrm{U}_{1}^{i-1}\mathrm{S}_{i+1}^{d}$ in Eq. \eqref{eq16} . However, for MHB, this number is $N_{\mathrm{MHB}}=n^{2}(\mathrm{S}_{1}^{d}\mathrm{U}_{1}^{d}+\mathrm{U}_{1}^{d}\mathrm{S}_{1}^{d}\mathrm{U}_{1}^{d})$.} To get an idea on the difference in computation consumption, supposing that $\mathrm{S}_{i}=2^5$ and $\mathrm{U}_{i}=11$, the ratio between the two numbers for comparison lists as
\textcolor{red}{\begin{equation}
	\frac{N_{\text{MHB}}}{N_{m\text{-} \mathrm{vCHB}}} = 
	\begin{cases}
		1 & d = 1 \\
		23.4 & d = 2 \\
		704.0 & d = 3 \\
		22203.7 & d = 4 \\
		\vdots & \vdots
	\end{cases},\label{eq39}
\end{equation}}
Here, the operation $\boldsymbol{\Gamma}_i{}^{0,-1}\partial_{\mathbf{Z}^{i-1}}\mathbf{F}^{i-1}\boldsymbol{\Gamma}_i{}_i^0$ with $\boldsymbol{\Gamma}_i{}^{0,-1}\in\mathbb{R}^{\mathrm{U}_i\times\mathrm{S}_i},\,\partial_{\mathbf{Z}^{i-1}}\mathbf{F}^{i-1}\in\mathbb{R}^{\mathrm{S}_i\times\mathrm{S}_i\times(n\mathrm{U}_1^{i-1})^2\mathrm{S}_{i+1}^d},\,\boldsymbol{\Gamma}_i{}_i^0\in\mathbb{R}^{\mathrm{S}_i\times\mathrm{U}_i}$ in Eq. \eqref{eq27} has the potential for parallelization, since the operation $\mathbb{R}^{\mathrm{U}_{i}\times\mathrm{S}_{i}}\times\mathbb{R}^{\mathrm{S}_{i}\times\mathrm{S}_{i}}\times\mathbb{R}^{\mathrm{S}_{i}\times U_{i}}$ is executed $(n\mathrm{U}_1^{i-1})^2\mathrm{S}_{i+1}^d$ times.

Thirdly, the robust phase condition is extended to fit the \textit{m}-VCF, which can be applied to fix the phase on the unknown base frequencies and be implemented by any one of the three methods, namely, HB, CO, and FD. In addition, the \textit{m}-VCF is also incorporated with two techniques, namely, the initialization of quasi-periodic solution after a NS bifurcation of periodic solutions in Section \ref*{sec4.2} and the stability analysis on quasi-periodic solutions by Lyapunov exponents in Appendix \ref{appD}. The flowchart of the \textit{m}-VCF is concluded in Appendix \ref{appE}.

%
%
%
%
%


\vspace{.5cm}
\noindent\textbf{Acknowledgments} Jun Jiang acknowledges the financial support of the National Natural Science Foundation of China (No. 12172267). Mingwu Li acknowledges the financial support of the National Natural Science Foundation of China (No. 12302014).

\newpage
\begin{appendices}

\section{Harmonic balance method}\label{appA}
In HB, the periodic functions $\boldsymbol{\Phi}_i(k_i,\tau_i)$ are defined by the Fourier series:
\begin{equation}
	\boldsymbol{\Phi}_i(k_i,\tau_i)=\begin{bmatrix}1,\cos(u_1\tau_i),\sin(u_1\tau_i),\cdots,\cos(u_{\mathrm{N}_i}\tau_i),\sin(u_{\mathrm{N}_i}\tau_i)\end{bmatrix}\in\mathbb{R}^{1\times\mathrm{U}_i},\label{EqA1}
\end{equation}
where $u_l,\,l=1,\cdots,\mathrm{N}_i$ is the harmonic order of the Fourier series and the term $1$ actually corresponds to the order of $u_0=0$. And the number of $\mathrm{U}_i$ is equal to  $2\mathrm{N}_i$. Here, the \textit{k}-parameter $k_{i;j}$ can be denoted as $u_0,\,u_l{}^c,\,u_l{}^s$, e.g., $\mathbf{Z}^{i}(u_{l}{}^{c},\boldsymbol{k}_{1}^{i-1},\boldsymbol{\tau}_{i+1}^{d})$ is the Fourier coefficients of $\mathbf{Z}^{i-1}(\boldsymbol{k}_{1}^{i-1},\boldsymbol{\tau}_{i}^{d})$ in term of the cosine term $\cos(u_l\tau_i)$. Notedly, the \textit{k}-parameter is chosen in a discretized frequency domain $k_i\in\mathbb{U}_i$. 

Actually, $\boldsymbol{\Phi}_i(k_i,\tau_i)$ and its derivatives can de expressed by:
\begin{equation}
	\begin{aligned}
		&\boldsymbol{\Phi}_{i}(k_{i},\tau_{i})=\boldsymbol{\Phi}_{i}(k_{i},\tau_{i})\boldsymbol{\nabla}_{i}{}^{0}
		\\&\boldsymbol{\Phi}_{i}^{\prime}(k_{i},\tau_{i})=\boldsymbol{\Phi}_{i}(k_{i},\tau_{i})\boldsymbol{\nabla}_{i}{}^{1}
		\\&\boldsymbol{\Phi}_{i}^{\prime\prime}(k_{i},\tau_{i})=\boldsymbol{\Phi}_{i}(k_{i},\tau_{i})\boldsymbol{\nabla}_{i}{}^{2}
	\end{aligned},\label{EqA2}
\end{equation}
where $\boldsymbol{\nabla}_{i}{}^{0},\,\boldsymbol{\nabla}_{i}{}^{1},\,\boldsymbol{\nabla}_{i}{}^{2}$ are three constant matrices if \textit{k}-parameter are fixed:
\begin{equation}
	\begin{aligned}
		&\boldsymbol{\nabla}_{i}{}_{i}^{0}=\mathbf{I}_{\mathrm{U}_{i}}
		\\&\boldsymbol{\nabla}_{i}{}^{1}=\mathrm{diagblk}\left(0,\boldsymbol{\nabla}_1,\cdots,\boldsymbol{\nabla}_{\mathrm{N}_i}\right),\boldsymbol{\nabla}_j=\begin{bmatrix}0&u_j\\-u_j&0\end{bmatrix}
		\\&\boldsymbol{\nabla}_{i}{}^{2}=\boldsymbol{\nabla}_{i}{}^1\boldsymbol{\nabla}_{i}{}^1
	\end{aligned},\label{EqA3}
\end{equation}
where the term diagblk means bock-diagonal matrix. 

The discretization of periodic functions $\boldsymbol{\Phi}_i(k_i,\tau_i)$ and derivatives are the transformation from the continuous hyper-time domain $\tau_{i}\in\mathbb{S}$ to the discretized frequency domain $k_i\in\mathbb{U}_i$. By using Fourier-Galerkin procedure
\begin{equation}
	\begin{aligned}
		&\frac{\left\lceil2\right\rceil}{2\pi}\int_{0}^{2\pi}\boldsymbol{\Phi}_{i}^{\mathrm{T}}\boldsymbol{\Phi}_{i}d\tau_{i}=\boldsymbol{\nabla}_{i}{}^{0}
		\\&\frac{\left\lceil2\right\rceil}{2\pi}\int_{0}^{2\pi}\boldsymbol{\Phi}_{i}^{\mathrm{T}}\boldsymbol{\Phi}_{i}^{\prime}d\tau_{i}=\boldsymbol{\nabla}_{i}{}^{1}
		\\&\frac{\left\lceil2\right\rceil}{2\pi}\int_{0}^{2\pi}\boldsymbol{\Phi}_{i}^{\mathrm{T}}\boldsymbol{\Phi}_{i}^{\prime\prime}d\tau_{i}=\boldsymbol{\nabla}_{i}{}^{2}
	\end{aligned},\label{EqA4}	
\end{equation}
one can obtain the results of Eq. \eqref{eq6}. Here $\left\lceil2\right\rceil$ is $2$ for $k_{i}\neq0$. And the matrices $\boldsymbol{\Upsilon}_i{}^0(k_i),\boldsymbol{\Upsilon}_i{}^1(k_i),\boldsymbol{\Upsilon}_i{}^2(k_i):\mathbb{U}_i\to\mathbb{R}^{\mathrm{U}_i\times\mathrm{U}_i}$ are $\boldsymbol{\nabla}_{i}{}^{0},\,\boldsymbol{\nabla}_{i}{}^{1},\,\boldsymbol{\nabla}_{i}{}^{2}\in\mathbb{R}^{\mathrm{U}_i\times\mathrm{U}_i}$.

Then, the transformation between the domain of $(\boldsymbol{k}_1^i,\overline{\boldsymbol{\tau}}_{i+1}^d)$ and $(\boldsymbol{k}_1^{i-1},\overline{\boldsymbol{\tau}}_{i}^d)$, the matrices $\boldsymbol{\Gamma}_i{}^0\in\mathbb{R}^{\mathrm{S}_i\times\mathrm{U}_i}$ and $\boldsymbol{\Gamma}_i{}^{0,1}\in\mathbb{R}^{\mathrm{U}_i\times\mathrm{S}_i}$ are matrices of inverse Discrete Fourier Transform (\textit{i}DFT) and Discrete Fourier Transform (DFT):
\begin{equation}
	\boldsymbol{\Gamma}_i{}^0=\boldsymbol{\Phi}_i(k_i,\overline{\boldsymbol{\tau}}_i)=\begin{bmatrix}1&\cos\left(u_1\tau_{i;1}\right)&\sin\left(u_1\tau_{i;1}\right)&\cdots&\cos\left(u_{\mathrm{N}_i}\tau_{i;1}\right)&\sin\left(u_{\mathrm{N}_i}\tau_{i;1}\right)\\\vdots&\vdots&\vdots&\ddots&\vdots&\vdots\\1&\cos\left(u_1\tau_{i;\mathrm{S}_i}\right)&\sin\left(u_1\tau_{i;\mathrm{S}_i}\right)&\cdots&\cos\left(u_{\mathrm{N}_i}\tau_{i;\mathrm{S}_i}\right)&\sin\left(u_{\mathrm{N}_i}\tau_{i;\mathrm{S}_i}\right)\end{bmatrix},\label{eqA5}
\end{equation}
\begin{equation}
	\left.\boldsymbol{\Gamma}_i{}^{0,-1}=\frac1{\mathrm{S}_i}\left[\begin{array}{ccc}1&\cdots&1\\2\cos\left(u_1\tau_{i;1}\right)&\cdots&2\cos\left(u_1\tau_{i;\mathrm{S}_i}\right)\\2\sin\left(u_1\tau_{i;1}\right)&\cdots&2\sin\left(u_1\tau_{i;\mathrm{S}_i}\right)\\\vdots&\ddots&\vdots\\2\cos\left(u_{\mathrm{N}_i}\tau_{i;1}\right)&\cdots&2\cos\left(u_{\mathrm{N}_i}\tau_{i;\mathrm{S}_i}\right)\\2\sin\left(u_{\mathrm{N}_i}\tau_{i;1}\right)&\cdots&2\sin\left(u_{\mathrm{N}_i}\tau_{i;\mathrm{S}_i}\right)\end{array}\right.\right],\label{eqA6}
\end{equation}
And the matrix $\boldsymbol{\Gamma}_i{}^1\in\mathbb{R}^{\mathrm{S}_i\times\mathrm{U}_i}$ is equal to $\boldsymbol{\Gamma}_i{}^1=\boldsymbol{\Gamma}_i{}^0\boldsymbol{\nabla}_{i}{}^{1}$. Last, the matrices $\boldsymbol{\Upsilon}_{\tau_i}{}^{1},\boldsymbol{\Upsilon}_{\tau_i}{}^{2}\in\mathbb{R}^{\mathrm{U}_i\times\mathrm{U}_i}$ in phase condition are $\boldsymbol{\nabla}_{i}{}^{1},\,\boldsymbol{\nabla}_{i}{}^{2}\in\mathbb{R}^{\mathrm{U}_i\times\mathrm{U}_i}$.

\section{Collocation method}\label{appB}
In CO, the domain of $k_i\in\mathbb{U}_i$ is same as the discretized time domain of $\overline{\boldsymbol{\tau}}_i\in\overline{\mathbb{S}}$ and the \textit{k}-parameter $k_{i;j}$ is defined as $u_{j}=2\left(j-1\right)\pi/\mathrm{U}_{i},\,j=1,2,\cdots,\mathrm{U}_{i}$, called as base points. Hence, the value of $\mathrm{U}_i$ is equal to $\mathrm{S}_i$ and $k_{i;j}=\tau_{i;j}$. For example, the coefficient $\mathbf{Z}^{i}(u_{j},\boldsymbol{k}_{1}^{i-1},\boldsymbol{\tau}_{i+1}^{d})$ is the value of $\mathbf{Z}^{i-1}(\boldsymbol{k}_{1}^{i-1},\tau_{i;j},\boldsymbol{\tau}_{i+1}^{d})$.

In CO, $\tau_i\in[0,2\pi)$ is divided into $\mathrm{P}_i$ equal-sized interval $\tau_{i,l}\in[2\pi m(l-1),2\pi ml)/\mathrm{U}_i$ with $m=\mathrm{U}_i/\mathrm{P}_i\in\mathbb{Z}_+$. The periodic functions $\boldsymbol{\Phi}_i(k_i,\tau_i)$ are also set as $\mathrm{P}_i$ set of piece-wise-functions $\boldsymbol{\Phi}_{i,l}(\boldsymbol{u}_{l},\tau_{i})$ by Lagrange basis polynomials functions
\begin{equation}
	\boldsymbol{\Phi}_i(k_i,\tau_i)=\boldsymbol{\Phi}_{i,l}(\boldsymbol{u}_l,\tau_i)=[0,\cdots0,L_{i,l,0},\cdots,L_{i,l,m},0,\cdots0,]\in\mathbb{R}^{1\times\mathrm{U}_i},\label{eqB1}
\end{equation}
when $\boldsymbol{u}_l=[u_{l,0},\cdots,u_{l,m}]\subset k_i$ with $u_{l,w}=u_{(l-1)m+w+1}=\tau_{i;(l-1)m+w+1},\,w=0,\cdots,m$. The function $L_{i,l,w}(\boldsymbol{u}_l,\tau_i)$ are the ($(l-1)m+w+1$)-th function in Eq. \eqref{eqB1}, except that $L_{i,\mathbf{P}_{i},0}(\boldsymbol{u}_{\mathbf{P}_{i}},\tau_{i})$ is the first function. And it is set as:
\begin{equation}
	L_{i,l,w}(\boldsymbol{u}_l,\tau_i)=
	\begin{cases}
		\prod_{o = 0, o\neq w}^{m} \frac{\tau_i - u_{l,o}}{u_{l,w} - u_{l,o}} & \tau_i \in \tau_{i,l}\\
		0 & \tau_i \notin \tau_{i,l}
	\end{cases}
	, \quad w = 0,\cdots,m,\label{eqB2}
\end{equation}
And the derivative of $L_{i,l,w}(\boldsymbol{u}_l,\tau_i)$ with respect to $\tau_i$ is:
\begin{equation}
	L'_{i,l,w}(\boldsymbol{u}_l,\tau_i)=
	\begin{cases}
		\sum_{o = 0, o\neq w}^{m} \left( \frac{1}{u_{l,w} - u_{l,o}} \prod_{p = 0, p\neq w, p\neq o}^{m} \frac{\tau_i - u_{l,p}}{u_{l,w} - u_{l,p}} \right) & \tau_i \in \tau_{i,l}\\
		0 & \tau_i \notin \tau_{i,l}
	\end{cases}
	, \quad w = 0,\cdots,m,\label{eqB3}
\end{equation}
The discretization of periodic functions $\boldsymbol{\Phi}_{i}(k_{i},\tau_{i})$ and derivatives are the transformation from the continuous hyper-time domain $\tau_{i}\in\mathbb{S}$ to the discretized hyper-time domain $\tilde{\tau}_{i}\in\tilde{\mathbb{S}}$ at the collocation nodes $\tilde{\tau}_{i}$ by the defined domain $k_{i}\in\mathbb{U}_{i}$.
\begin{equation}
	\boldsymbol{\Phi}_{i}(k_{i},\tilde{\tau}_{i}) = \tilde{\boldsymbol{\mathcal{L}}}_{i}{}^{0} = 
	\begin{bmatrix}
		\tilde{L}_{i,1,0} & \cdots & \tilde{L}_{i,1,m} \\
		& &\tilde{L}_{i,2,0} & \cdots & \tilde{L}_{i,2,m} \\
		& & & & \tilde{L}_{i,3,0} & \cdots & \tilde{L}_{i,3,m} \\
		& & & & & &\ddots \\
		\tilde{L}_{i,P_i,m} & & & & & & \tilde{L}_{i,P_i,0} & \cdots 
	\end{bmatrix},\label{eqB4}
\end{equation}
where $\tilde{L}_{i,l,w} = L_{i,l,w}(\boldsymbol{u}_l,\tilde{\tau}_{i,l}) \in \mathbb{R}^{m \times 1}$, with $\tilde{\tau}_{i,l}$ are the $m$ collocation nodes in $\tau_{i,l}$. $\boldsymbol{\Phi}'_{i}(k_{i},\tilde{\tau}_{i}) = \tilde{\boldsymbol{\mathcal{L}}}_{i}{}^{1} \in \mathbb{R}^{\mathrm{U_i} \times \mathrm{U_i}}$ has the same construction as $\tilde{\boldsymbol{\mathcal{L}}}_{i}{}^{0} \in \mathbb{R}^{\mathrm{U_i} \times \mathrm{U_i}}$, and $\boldsymbol{\Phi}''_{i}(k_{i},\tilde{\tau}_{i}) = \tilde{\boldsymbol{\mathcal{L}}}_{i}{}^{1} \overline{\boldsymbol{\mathcal{L}}}_{i}{}^{1} \in \mathbb{R}^{\mathrm{U_i} \times \mathrm{U_i}}$, based on the relationship of $\boldsymbol{\Phi}''_{i}(k_{i},\tilde{\tau}_{i}) = \boldsymbol{\Phi}'_{i}(k_{i},\tilde{\tau}_{i}) \boldsymbol{\Phi}'_{i}(k_{i},\overline{\tau}_{i})$. And the matrices $\boldsymbol{\Upsilon}_{i}{}^{0}(k_{i}), \boldsymbol{\Upsilon}_{i}{}^{1}(k_{i}), \boldsymbol{\Upsilon}_{i}{}^{2}(k_{i}): U_i \to \mathbb{R}^{\mathrm{U_i} \times \mathrm{U_i}}$ are $\tilde{\boldsymbol{\mathcal{L}}}_{i}{}^{0}, \tilde{\boldsymbol{\mathcal{L}}}_{i}{}^{1}, \tilde{\boldsymbol{\mathcal{L}}}_{i}{}^{2} \in \mathbb{R}^{\mathrm{U_i} \times \mathrm{U_i}}$.

Then, the transformation between the domain of $(\boldsymbol{k}^{i}_{1},\overline{\boldsymbol{\tau}}^{d}_{i + 1})$ and $(\boldsymbol{k}^{i - 1}_{1},\overline{\boldsymbol{\tau}}^{d}_{i})$, the matrices $\boldsymbol{\Gamma}_i{}^0\in\mathbb{R}^{\mathrm{S}_i\times\mathrm{U}_i}$ and $\boldsymbol{\Gamma}_i{}^{0,1}\in\mathbb{R}^{\mathrm{U}_i\times\mathrm{S}_i}$ are identity matrices with size of $\mathrm{U}_i\times\mathrm{S}_i$. And the matrix $\boldsymbol{\Gamma}_i{}^1\in\mathbb{R}^{\mathrm{S}_i\times\mathrm{U}_i}$ is equal to $\overline{\boldsymbol{\mathcal{L}}}_{i}{}^{1}$. Last, the matrices $\boldsymbol{\Upsilon}_{\tau_i}{}^{1},\boldsymbol{\Upsilon}_{\tau_i}{}^{2}\in\mathbb{R}^{\mathrm{U}_i\times\mathrm{U}_i}$ in phase condition are $\overline{\boldsymbol{\mathcal{L}}}_{i}{}^{1},\,\overline{\boldsymbol{\mathcal{L}}}_{i}{}^{2}\in\mathbb{R}^{\mathrm{U}_i\times\mathrm{U}_i}$ with $\overline{\boldsymbol{\mathcal{L}}}_{i}{}^{2}=\overline{\boldsymbol{\mathcal{L}}}_{i}{}^{1}\overline{\boldsymbol{\mathcal{L}}}_{i}{}^{1}$

\section{Finite difference method}\label{appC}

In FD, the domain of $k_{i} \in \mathbb{U}_{i}$ is also same as the discretized time domain of $\overline{\tau}_{i} \in \overline{\mathbb{S}}$. But there are no expressions of periodic functions $\boldsymbol{\Phi}_{i}(k_{i},\tau_{i})$. The discretization is the process from $\mathbf{Z}^{i - 1}(\boldsymbol{\tau}_{i}^{d}), \dot{\mathbf{Z}}^{i - 1}(\boldsymbol{\tau}_{i}^{d}), \ddot{\mathbf{Z}}^{i - 1}(\boldsymbol{\tau}_{i}^{d})$ to $\mathbf{Z}^{i - 1}(\overline{\boldsymbol{\tau}}_{i},\boldsymbol{\tau}_{i + 1}^{d}), \dot{\mathbf{Z}}^{i - 1}(\overline{\boldsymbol{\tau}}_{i},\boldsymbol{\tau}_{i}^{d}), \ddot{\mathbf{Z}}^{i - 1}(\overline{\boldsymbol{\tau}}_{i},\boldsymbol{\tau}_{i}^{d})$.
\begin{equation}
	\begin{aligned}
		&\mathbf{Z}^{i-1}(\overline{\boldsymbol{\tau}}_{i},\boldsymbol{\tau}_{i}^{d})
		\approx
		[\mathbf{I}_{n\mathrm{U}_{1}^{i-1}}\otimes\overline{\boldsymbol{\mathcal{T}}}_{i}{}^{0}(k_{i})]
		\mathbf{Z}^{i}(k_{i},\boldsymbol{\tau}_{i+1}^{d})
		\\&\dot{\mathbf{Z}}^{i-1}(\overline{\boldsymbol{\tau}}_{i},\boldsymbol{\tau}_{i}^{d})
		\approx
		[\mathbf{I}_{n\mathrm{U}_{1}^{i-1}}\otimes\omega_{i}\overline{\boldsymbol{\mathcal{T}}}_{i}{}^{1}(k_{i})]
		\mathbf{Z}^{i}(k_{i},\tau_{i+1}^{d})
		+[\mathbf{I}_{n\mathrm{U}_{1}^{i-1}}\otimes\overline{\boldsymbol{\mathcal{T}}}_{i}{}^{0}(k_{i})]
		\dot{\mathbf{Z}}^i(k_i,\boldsymbol{\tau}_{i+1}^d)
		\\&\ddot{\mathbf{Z}}^{i-1}(\overline{\boldsymbol{\tau}}_{i},\boldsymbol{\tau}_{i}^{d})
		\approx
		[\mathbf{I}_{n\mathrm{U}_{1}^{i-1}}\otimes\omega_{i}{}^2\overline{\boldsymbol{\mathcal{T}}}_{i}{}^{2}(k_{i})]
		\mathbf{Z}^{i}(k_{i},\tau_{i+1}^{d})
		+2[\mathbf{I}_{n\mathrm{U}_{1}^{i-1}}\otimes\omega_{i}\overline{\boldsymbol{\mathcal{T}}}_{i}{}^{1}(k_{i})]
		\dot{\mathbf{Z}}^i(k_i,\tau_{i+1}^d)
		\\&+[\mathbf{I}_{n\mathrm{U}_{1}^{i-1}}\otimes\overline{\boldsymbol{\mathcal{T}}}_{i}{}^{0}(k_{i})]
		\ddot{\mathbf{Z}}^i(k_i,\boldsymbol{\tau}_{i+1}^d)
	\end{aligned},\label{eqC1}
\end{equation}
where $\overline{\boldsymbol{\mathcal{T}}}_{i}{}^{0}, \overline{\boldsymbol{\mathcal{T}}}_{i}{}^{1}, \overline{\boldsymbol{\mathcal{T}}}_{i}{}^{2} \in \mathbb{R}^{\mathrm{U_i} \times \mathrm{U_i}}$ are the discretized results of $\boldsymbol{\Phi}_i(k_i,\tau_i)$ and its derivatives at $\tau_i=\overline{\boldsymbol{\tau}}_i\in\overline{\mathbb{S}}$. Here, $\overline{\boldsymbol{\mathcal{T}}}_{i}{}^{0}$ is an identity matrix and $\overline{\boldsymbol{\mathcal{T}}}_{i}{}^{1},\,\overline{\boldsymbol{\mathcal{T}}}_{i}{}^{2}$ are constructed by:
\begin{equation}
	\overline{\boldsymbol{\mathcal{T}}}_{i}{}^{g}=\frac1{(\Delta\tau_i)^\mathrm{g}}\begin{bmatrix}\alpha_{_{\mathrm{g},0}}&\cdots&\alpha_{_{\mathrm{g},+r}}&&&\alpha_{_{\mathrm{g},-l}}&\cdots\\\cdots&\alpha_{_{\mathrm{g},0}}&\cdots&\ddots&&&\alpha_{_{\mathrm{g}}}\\\alpha_{_{\mathrm{g},-l}}&\cdots&\alpha_{_{\mathrm{g},0}}&\ddots&\alpha_{_{\mathrm{g},+r}}&&\\&\alpha_{_{\mathrm{g},-l}}&\cdots&\ddots&\cdots&\alpha_{_{\mathrm{g},+r}}&\\&&\alpha_{_{\mathrm{g},-l}}&\ddots&\alpha_{_{\mathrm{g},0}}&\cdots&\alpha_{_{\mathrm{g},+r}}\\\alpha_{_{\mathrm{g},+r}}&&&\ddots&\cdots&\alpha_{_{\mathrm{g},0}}&\cdots\\\cdots&\alpha_{_{\mathrm{g},+r}}&&&\alpha_{_{\mathrm{g},-l}}&\cdots&\alpha_{_{\mathrm{g},0}}\end{bmatrix},\,g=1,2,\label{eqC2}
\end{equation}
where the values of $\alpha_{g,0}, \alpha_{g, + r}, \alpha_{g, - l}$ can be determined by the following procedures. First, the derivatives of $\boldsymbol{Z}^{i-1}(\boldsymbol{\tau}_{i}^{d})$ are asked to equal to:
\begin{equation}
	\left.\frac{\partial^g\mathbf{Z}^{i-1}(\tau_i,\boldsymbol{\tau}_{i+1}^d)}{\partial\tau_i{}^g}\right|_{\tau_{i;j}}=\frac{1}{(\Delta\tau_i)^g}\sum_{p=-l}^{+r}\alpha_{g,p}\mathbf{Z}^i(k_{i;j+p},\boldsymbol{\tau}_{i+1}^d),\label{eqC3}
\end{equation}
based on the arbitrary number of adjacent base points $K = [-l,\cdots,0,\cdots,+r] \in \mathbb{R}^{1 \times N_K}$ with $r, l \geq 0$, $N_K \leq U_i$, where $\boldsymbol{Z}^{i}(k_{;j + p},\boldsymbol{\tau}_{i + 1}^{d}) = \boldsymbol{Z}^{i - 1}(\tau_{i;j + p},\boldsymbol{\tau}_{i + 1}^{d})$ and it can be expanded by the Taylor expression:
\begin{equation}
	\mathbf{Z}^i(k_{i;j+p},\boldsymbol{\tau}_{i+1}^d)=\sum_{m=0}^\infty\frac{1}{m!}\left.\frac{\partial^m\mathbf{Z}^{i-1}(\tau_i^d)}{\partial\tau_i{}^m}\right|_{\tau_{i;j}}(p\Delta\tau_i)^m,\label{eqC4}
\end{equation}
So $\alpha_{g,p}$ can be computed by:
\begin{equation}
	\begin{bmatrix}0\\\vdots\\\mathrm{g}!\\\vdots\\0\end{bmatrix}=\begin{bmatrix}1&\cdots&1&\cdots&1\\\vdots&\ddots&\vdots&\ddots&\vdots\\\left(-l\right)^\mathrm{g}&\cdots&0^\mathrm{g}&\cdots&r^\mathrm{g}\\\vdots&\ddots&\vdots&\ddots&\vdots\\\left(-l\right)^{N_K}&\cdots&0^{N_K}&\cdots&r^{N_K}\end{bmatrix}\begin{bmatrix}\alpha_{\mathrm{g},-l}\\\vdots\\\alpha_{\mathrm{g},0}\\\vdots\\\alpha_{\mathrm{g},+r}\end{bmatrix},\label{eqC5}
\end{equation}
And the matrices $\boldsymbol{\Upsilon}_{i}{}^{0}(k_{i}), \boldsymbol{\Upsilon}_{i}{}^{1}(k_{i}), \boldsymbol{\Upsilon}_{i}{}^{2}(k_{i}): U_i \to \mathbb{R}^{\mathrm{U_i} \times \mathrm{U_i}}$ are $\overline{\boldsymbol{\mathcal{T}}}_{i}{}^{0}, \overline{\boldsymbol{\mathcal{T}}}_{i}{}^{1}, \overline{\boldsymbol{\mathcal{T}}}_{i}{}^{2} \in \mathbb{R}^{\mathrm{U_i} \times \mathrm{U_i}}$. The matrices $\boldsymbol{\Gamma}_i{}^0\in\mathbb{R}^{\mathrm{S}_i\times\mathrm{U}_i}$ and $\boldsymbol{\Gamma}_i{}^{0,1}\in\mathbb{R}^{\mathrm{U}_i\times\mathrm{S}_i}$ are identity matrices with size of $\mathrm{U}_i\times\mathrm{S}_i$. And the matrix $\boldsymbol{\Gamma}_i{}^1\in\mathbb{R}^{\mathrm{S}_i\times\mathrm{U}_i}$ is equal to $\overline{\boldsymbol{\mathcal{T}}}_{i}{}^{1}$. Last, the matrices $\boldsymbol{\Upsilon}_{\tau_i}{}^{1},\boldsymbol{\Upsilon}_{\tau_i}{}^{2}\in\mathbb{R}^{\mathrm{U}_i\times\mathrm{U}_i}$ in phase condition are $\overline{\boldsymbol{\mathcal{T}}}_{i}{}^{1},\,\overline{\boldsymbol{\mathcal{T}}}_{i}{}^{2}\in\mathbb{R}^{\mathrm{U}_i\times\mathrm{U}_i}$.

\section{Lyapunov exponents}\label{appD}

Assuming that $\mathbf{x}_s\left(t\right)$ is a stationary solution of Eq. (32), $\dot{\Psi}(t)\approx\mathbf{J}(t)\Psi(t),\,\mathbf{\Psi}(0)=\mathbf{I}_{2n}$ is the perturbation system which reveals the evolution of perturbation $\mathbf{\Psi}(0)$ over time. Here, $\mathbf{J}\left(t\right)$ is a Jacobian matrix, governed by the stationary solution: 
\begin{equation}
	\mathbf{J}\left(t\right)=
	\begin{bmatrix}\mathbf{0}
		&\mathbf{I}_n\\-\mathbf{M}^{0,-1}\left(\mathbf{K}^0+\frac{\partial\mathbf{F}^0}{\partial\mathbf{Z}^0}\bigg|_{\mathbf{x}_s(t)}\right)&-\mathbf{M}^{0,-1}\left(\mathbf{D}^0+\frac{\partial\mathbf{F}^0}{\partial\dot{\mathbf{Z}^0}}\bigg|_{\mathbf{x}_s(t)}\right)
	\end{bmatrix}_{2n\times2n},\label{eqD1}
\end{equation}
The dimension \textit{d} of the base frequencies $\boldsymbol{\omega}_1^d$ can be used to categorize the methods to assess stability:

1. Equilibrium solution: $\mathbf{x}_s\left(t\right)\equiv\mathbf{x}_s\left(0\right)$:

In this case, the Jacobian matrix is also constant. The evolution of perturbation system is actually governed by the Jacobian matrix. The stability of equilibrium point is evaluated by the eigenvalues of Jacobian matrix \cite{JW41}, where if the real parts of all eigenvalues is not greater than 0, it is stable; otherwise, it is unstable.

2. Periodic solution: $\mathbf{x}_s\left(\omega,t+T\right)\equiv\mathbf{x}_s\left(\omega,t\right), T=\frac{2\pi}\omega $:

Due to the periodicity of periodic solution, the Jacobian matrix also has periodicity $\mathbf{J}\left(t+T\right)=\mathbf{J}\left(t\right)$. The map of perturbation is expressed by $\Delta\mathbf{x}\left(kT\right)=\mathbf{\Psi}\left(T,0\right)^k\Delta\mathbf{x}\left(0\right)$, where the constant matrix $\mathbf{M}=\mathbf{\Psi}(T,0)$ is called as monodromy matrix and is generally used to evaluate stability of periodic response (Floquet theory) \cite{JW42}. If the absolute of all eigenvalues (also called as Floquet multipliers) of monodromy matrix is not larger than 1, the periodic response is stable; otherwise, it is unstable. The Hill's method \cite{JW43,JW44,JW45} are also can be used to evaluate the stability of periodic solutions, where the stabilities of the Fourier coefficient in HB are considered as the stabilities of periodic solutions.

3. Quasi-periodic solution ($d\geq2$): $\mathbf{x}_s\left(t\right)=\mathbf{x}_s\left(\boldsymbol{\omega}_1^d,t\right)$:

Since the quasi-periodic solution $\mathbf{x}_s\left(t\right)$ never repeats in the time domain, the Jacobian matrix also never repeats. The Lyapunov exponent based on the discrete Gram-Schmidt orthonormalization \cite{JW46,JW47} is used to approximately assess the stability of a quasi-periodic solution:
\begin{equation}
	\begin{aligned}&\sigma^{(m)}\left(i\Delta t\right)=\lim_{i\to\infty}\frac1{i\Delta t}\ln\mathrm{V}_m\left(i\Delta t\right), m=1,2,\cdots,2n\\&\sigma_m\left(i\Delta t\right)=\sigma^{(m)}\left(i\Delta t\right)-\sigma^{(m-1)}\left(i\Delta t\right)\end{aligned},\label{eqD2}
\end{equation}
where $\sigma^{(m)}(i\Delta t)$ is the \textit{m}-th order Lyapunov exponent, $\sigma_m(i\Delta t)$ is the \textit{m}-th first-order Lyapunov exponent. If all first-order Lyapunov exponents is not larger than 0, the quasi-periodic solution is stable; otherwise, it is unstable.

Notedly, $\mathrm{V}_m\left(i\Delta t\right)$ is the volume of \textit{m}-dimensional parallelepiped of the transition matrix $\mathbf{\Psi}(i\Delta t,0)$. Here, $\mathbf{\Psi}(i\Delta t,0)$ is calculated by $\mathbf{\Psi}\left(i\Delta t,\left(i-1\right)\Delta t\right)\hat{\mathbf{\Psi}}\left(\left(i-1\right)\Delta t,0\right)$, where $\hat{\mathbf{\Psi}}\left(\left(i-1\right)\Delta t,0\right)$ is re-orthonormalized and re-normalized by the transition matrix $\mathbf{\Psi}\left(\left(i-1\right)\Delta t,0\right)$ based on the discrete Gram-Schmidt orthonormalization.

Given the quasi-periodicity of $\mathbf{x}_s\left(\boldsymbol{\omega}_1^d,t\right)$, defining the time interval $\Delta t$ as $\Delta t=T_j=\frac{2\pi}{\omega_j},j=1,2,\cdots,d$, the transition matrix is transformed into $\mathbf{\Psi}\left(iT_j,(i-1)T_j\right):=\mathbf{\Psi}\left(\boldsymbol{\tau}_{e;i},\boldsymbol{\tau}_{s;i}\right)$ form the start point $\boldsymbol{\tau}_{s;i}$ to the end point $\boldsymbol{\tau}_{e;i}$, where these two points are set as:
\begin{equation}
	\begin{aligned}&\boldsymbol{\tau}_{s;i}=\left[\left(i-1\right)\boldsymbol{\omega}_1^{j-1}T_j\;mod\;2\pi;0;\left(i-1\right)\boldsymbol{\omega}_{j+1}^dT_j\;mod\;2\pi\right]\\&\boldsymbol{\tau}_{e;i}=\left[i\boldsymbol{\omega}_1^{j-1}T_j\;mod\;2\pi;2\pi;i\boldsymbol{\omega}_{j+1}^dT_j\;mod\;2\pi\right]\end{aligned},\label{eqD3}
\end{equation}
Here, $\boldsymbol{\tau}_{s;i}$ fills a subset of the hyper-time domain $\mathbb{T}^d$ densely for $i\to+\infty$, denoted as:
\begin{equation}
	\mathbb{S}_s:=\left[\tau_1,\cdots,\tau_{j-1},\tau_{j+1},\cdots,\tau_d\right],\label{eqD4}
\end{equation}
Fiedler and Hetzler \cite{JW48} recently proposed an efficient numerical calculation for approximating this transition matrix, where $\mathrm{Y}=\prod_{k=1,k\neq j}^d\mathrm{Y}_k$ points are sampled in the subset $\mathbb{S}_s$. These \textit{d}-1 values of $\mathrm{Y}_k,\,k=1,\cdots,j-1,j+1,\cdots,d$ can be recognized as \textit{d}-1 values of $\mathrm{S}_{i},\,i=1,2,\cdots,d-1$ to get a set of $\mathrm{Y}$ points, denoted as $\bar{\boldsymbol{\tau}}_1^{d-1}$ in Eq. \eqref{eq17a}. Choosing the q-th point of $\bar{\boldsymbol{\tau}}_1^{d-1}$ defines a start point $\boldsymbol{\tau}_{s,q}:=\left[\tau_{s,q,1};\cdots;\tau_{s,q,j-1};0;\tau_{s,q,j+1},\cdots,\tau_{s,q,d}\right],\,q=1,2,\cdots,\mathrm{Y}$, the corresponding end point is $\boldsymbol{\tau}_{s,e}$, where $\tau_{e,q,k}=\tau_{s,q,k}+\omega_kT_j\,mod\,2\pi,\,k\neq j$ and $\tau_{e,q,j}=2\pi$. Using the method of matrix exponential approximation, the transition matrix $\mathbf{\Psi}\left(\boldsymbol{\tau}_{e,q},\boldsymbol{\tau}_{s,q}\right)$ is calculated by:
\begin{equation}
	\mathbf{\Psi}\left(\boldsymbol{\tau}_{e,q},\boldsymbol{\tau}_{s,q}\right)\approx\prod_{\eta=1}^{N_{\mathrm{M}}}\mathrm{e}^{[\mathbf{J}_{\eta+1/2}]\frac{T_{j}}{N_{\mathrm{M}}}},\mathbf{J}_{\eta+1/2}=\frac{\mathbf{J}(\boldsymbol{\tau}_{\eta+1})+\mathbf{J}(\boldsymbol{\tau}_{\eta})}{2},\label{eqD5}
\end{equation}
where $N_{\mathrm{M}}$ is number of sub-intervals in $T_{j}$, and the \textit{k}-th element of $\boldsymbol{\tau}_{\eta}$ is the $\eta$-th point defined as:
\begin{equation}
	\tau_{\eta,k}=\tau_{s,q,k}+\omega_k\frac{\eta-1}{N_\mathrm{M}}T_j\,mod\,2\pi,\, k=1,2,\cdots,d,\label{eqD6}
\end{equation}
As a consequence, the quasi-periodic matrix $\mathbf{\Psi}\left(\boldsymbol{\tau}_{s}\right)$ can be expressed by means of trigonometric polynomials based on (\textit{d}-1)-dimensional FFT. In this case, the transition matrix $\mathbf{\Psi}\left(iT_j,(i-1)T_j\right):=\mathbf{\Psi}\left(\boldsymbol{\tau}_{e;i},\boldsymbol{\tau}_{s;i}\right)$ can be approximated by the function $\mathbf{\Psi}\left(\boldsymbol{\tau}_{s}\right)$ with $\boldsymbol{\tau}_s=\boldsymbol{\tau}_{s;i}$. In this work, the maximum of \textit{i} is denoted as $N_{\mathrm{L}}$.

\section{Flowchart of the \textit{m}-VCF}\label{appE}
\begin{figure*}[htbp]
	\centering
	\includegraphics[width=0.9\textwidth]{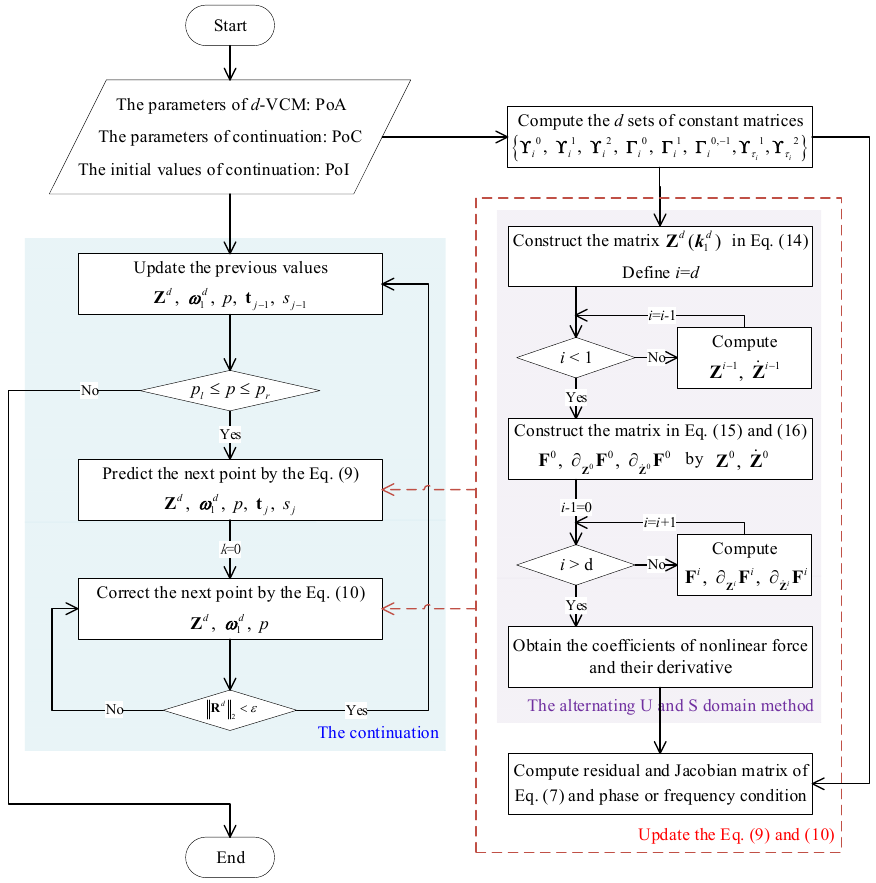}
	\caption{The flowchart of the \textit{m}-VCF.}\label{figE1}
\end{figure*}

In this appendix, the flowchart of the m-VCF is shown in Fig. \ref*{figE1}. 

\end{appendices}


\begin{thebibliography}{99}
	
	\bibitem {JW1}P. Baldi, R. Montalto, Quasi-periodic incompressible Euler flows in 3D, \href{https://doi.org/10.1016/j.aim.2021.107730}{Adv. Math.} 384, 107730 (2021).
	
	\bibitem {JW2}A. Cheffer, M.A. Savi, Biochaos in cardiac rhythms, \href{https://doi.org/10.1140/epjs/s11734-021-00314-7}{Eur. Phys. J. Spec. Top.} 231, 833-845 (2022).
	
	\bibitem {JW3}S. Fan, L. Hong, J. Jiang, Blue-Sky Catastrophic Bifurcations Behind Emergence and Disappearance of Quasiperiodic Rubbing Oscillations in a Piecewise Smooth Rotor–Stator System, \href{https://doi.org/10.1142/s0218127422502212}{Int. J. Bifurcation Chaos}, 32 (2022) 2250221.
	
	\bibitem {JW4}A. Mereles, D.S. Alves, K.L. Cavalca, On the continuation of quasi-periodic solutions of rotor systems with fluid-film bearings, \href{https://doi.org/10.1007/s11071-024-10616-9}{Nonlinear Dyn.}, 113 (2025) 9645–9665.
	
	\bibitem {JW5}T.M. Wu, J.L. Huang, W.D. Zhu, Quasi-periodic oscillation characteristics of a nonlinear energy sink system under harmonic excitation, \href{https://doi.org/10.1016/j.jsv.2023.118143}{J.Sound Vib.}, 572 (2024) 118143.
	
	\bibitem {JW6}N. Coudeyras, S. Nacivet, J.J. Sinou, Periodic and quasi-periodic solutions for multi-instabilities involved in brake squeal, \href{https://doi.org/10.1016/j.jsv.2009.08.017}{J.Sound Vib.}, 328 (2009) 520-540.
	
	\bibitem {JW49}L. Lan, Y. Ni, Y. Jiang and J. Li, Motion of the Moonlet in the Binary System 243 Ida, \href{https://doi.org/10.1007/s10409-017-0722-3}{Acta Mech. Sin.}, 34 (2018) 224.
	
	\bibitem {JW7}M. Guskov, F. Thouverez, Harmonic balance-based approach for quasi-periodic motions and stability analysis, \href{https://doi.org/10.1115/1.4005823}{J. Vib. Acoust.}, 134 (2012) 031003.
	
	\bibitem {JW8}K. Prabith, I.R.P. Krishna, Bifurcation studies of a nonlinear mechanical system subjected to multi-frequency-quasi-periodic excitations, in: Advances in Nonlinear Dynamics, (Springer, Cham, 2020), \href{https://doi.org/10.1007/978-3-030-81162-4_63}{https://doi.org/10.1007/978-3-030-81162-4-63}.
	
	\bibitem {JW9}K. Prabith, I.R.P. Krishna, A time variational method for the approximate solution of nonlinear systems undergoing multiple-frequency excitations, \href{https://doi.org/10.1115/1.4045944}{J. Comput. Nonlinear Dynam.}, 15 (2020) 031006.
	
	\bibitem {JW10}C. Kaas-Petersen, Computation of quasi-periodic solutions of forced dissipative systems II, \href{https://doi.org/10.1016/0021-9991(86)90042-2}{J. Comput. Phys.}, 64 (1986) 433-442.
	
	\bibitem {JW11}C. Kaas-Petersen, Computation, continuation, and bifurcation of torus solutions for dissipative maps and ordinary differential equations, \href{https://doi.org/10.1016/0167-2789(87)90105-9}{Physica D}, 25 (1987) 288-306.
	
	\bibitem {JW12}Y.B. Kim, Quasi-periodic response and stability analysis for non-linear systems: a general approach, \href{https://doi.org/10.1006/jsvi.1996.0220}{J.Sound Vib.}, 192 (1996) 821-833.
	
	\bibitem {JW13}Y.B. Kim, S.T. Noah, Quasi-periodic response and stability analysis for a non-linear jeffcott rotor, \href{https://doi.org/10.1006/jsvi.1996.0059}{J.Sound Vib.}, 190 (1996) 239-253.
	
	\bibitem {JW14}L. Junge, C. Frey, G. Ashcroft, and E. Kügeler, A new Harmonic Balance approach using multidimensional time, \href{https://doi.org/10.1115/1.4049698}{J. Eng. Gas Turbines Power.}, 143 (2021).
	
	\bibitem {JW15}T.M. Cameron, J.H. Griffin, An alternating frequency/time domain method for calculating the steady-state response of nonlinear dynamic systems, \href{https://doi.org/10.1115/1.3176036}{J. Appl. Mech.}, 56 (1989) 149-154.
	
	\bibitem {JW16}F. Schilder, H.M. Osinga, W. Vogt, Continuation of Quasi-periodic Invariant Tori, \href{https://doi.org/10.1137/040611240}{SIAM J. Appl. Dyn. Syst.}, 4 (2005) 459-488.
	
	\bibitem {JW17}À. Haro, A. Luque, J.M. Mondelo, M. Canadell, and J.L. Figueras, The parameterization method for invariant manifolds (Springer, Cham, 2016), \href{https://doi.org/10.1007/978-3-319-29662-3}{https://doi.org/10.1007/978-3-319-29662-3}.
	
	\bibitem {JW18}S. Bäuerle, R. Fiedler, H. Hetzler, An engineering perspective on the numerics of quasi-periodic oscillations, \href{https://doi.org/10.1007/s11071-022-07407-5}{Nonlinear Dyn.}, 108 (2022) 3927-3950.
	
	\bibitem {JW19}H. Liao, Q. Zhao, D. Fang, The continuation and stability analysis methods for quasi-periodic solutions of nonlinear systems, \href{https://doi.org/10.1007/s11071-020-05497-7}{Nonlinear Dyn.}, 100 (2020) 1469-1496.
	
	\bibitem {JW20}J. Wu, L. Hong, J. Jiang, A comparative study on multi- and variable-coefficient harmonic balance methods for quasi-periodic solutions, \href{https://doi.org/10.1016/j.ymssp.2022.109929}{Mech. Syst. Sig. Process.}, 187 (2023) 109929.
	
	\bibitem {JW21}N.N. Balaji, J. Gross, M. Krack, Harmonic Balance for quasi-periodic vibrations under nonlinear hysteresis, \href{https://doi.org/10.1016/j.jsv.2024.118570}{J.Sound Vib.}, 590 (2024) 118570.
	
	\bibitem {JW22}G. Ghannad Tehrani, C. Gastaldi, T.M. Berruti, Trained Harmonic Balance Method for Parametrically Excited Jeffcott Rotor Analysis, \href{https://doi.org/10.1115/1.4048578}{J. Comput. Nonlinear Dynam.}, 16 (2020) 011003.
	
	\bibitem {JW23}L. Guillot, P. Vigué, C. Vergez, B. Cochelin, Continuation of quasi-periodic solutions with two-frequency Harmonic Balance Method, \href{https://doi.org/10.1016/j.jsv.2016.12.013}{J.Sound Vib.}, 394 (2017) 434-450.
	
	\bibitem {JW24}J.L. Huang, W.D. Zhu, An incremental harmonic balance method with two timescales for quasiperiodic motion of nonlinear systems whose spectrum contains uniformly spaced sideband frequencies, \href{https://doi.org/10.1007/s11071-017-3708-6}{Nonlinear Dyn.}, 90 (2017) 1015-1033.
	
	\bibitem {JW25}J.L. Huang, W.D. Zhu, A New Incremental Harmonic Balance Method With Two Time Scales for Quasi-Periodic Motions of an Axially Moving Beam With Internal Resonance Under Single-Tone External Excitation, \href{https://doi.org/10.1115/1.4035135}{J. Vib. Acoust.}, 139 (2017).
	
	\bibitem {JW26}B. Zhou, F. Thouverez, D. Lenoir, A variable-coefficient harmonic balance method for the prediction of quasi-periodic response in nonlinear systems, \href{https://doi.org/10.1016/j.ymssp.2015.04.022}{Mech. Syst. Sig. Process.}, 64-65 (2015) 233-244.
	
	
	\bibitem {JW27}Z. Zheng, Z. Lu, G. Liu, Y. Chen, Twice Harmonic Balance Method for Stability and Bifurcation Analysis of Quasi-Periodic Responses, \href{https://doi.org/10.1115/1.4055923}{J. Comput. Nonlinear Dynam.}, 17 (2022) 121006.
	
	
	\bibitem {JW28}H. Hetzler, S. Bäuerle, Stationary solutions in applied dynamics: A unified framework for the numerical calculation and stability assessment of periodic and quasi-periodic solutions based on invariant manifolds, \href{https://doi.org/10.1002/gamm.202300006}{GAMM-Mitteilungen}, n/a (2023) e202300006.
	
	
	\bibitem {JW29}D. Roose, R. Szalai, Continuation and Bifurcation Analysis of Delay Differential Equations, in: Numerical Continuation Methods for Dynamical Systems: Path following and boundary value problems, (Springer, Dordrecht, 2007), \href{https://doi.org/10.1007/978-1-4020-6356-5_12}{https://doi.org/10.1007/978-1-4020-6356-5-12}.
	
	
	\bibitem {JW30}H. Dankowicz, F. Schilder, Recipes for Continuation, (SIAM, Philadelphia, 2013), \href{https://doi.org/10.1137/1.9781611972573}{https://doi.org/10.1137/1.9781611972573}.
	
	\bibitem {JW31}M. Li, Tor: a toolbox for the continuation of two-dimensional tori in autonomous systems and non-autonomous systems with periodic forcing, (2020), \href{https://10.48550/arXiv.2012.13256}{https://10.48550/arXiv.2012.13256}.
	
	\bibitem {JW32}M. Krack, J. Gross, Harmonic balance for nonlinear vibration problems, (Springer, Cham, 2019), \href{https://doi.org/10.1007/978-3-030-14023-6}{https://doi.org/10.1007/978-3-030-14023-6}.
	
	\bibitem {JW33}Z. Yan, H. Dai, Q. Wang, S.N. Atluri, Harmonic balance methods: A review and recent developments, \href{https://doi.org/10.32604/cmes.2023.028198}{CMES-Comp. Model. Eng. Sci.}, 137 (2023) 1419-1459.
	
	\bibitem {JW34}R.J. Kuether, A. Steyer, Large-scale harmonic balance simulations with Krylov subspace and preconditioner recycling, \href{https://doi.org/10.1007/s11071-023-09171-6}{Nonlinear Dyn.}, 112 (2024) 3377-3398.
	
	\bibitem {JW35}T. Vadcard, F. Thouverez, A. Batailly, On the detection of nonlinear normal mode-related isolated branches of periodic solutions for high-dimensional nonlinear mechanical systems with frictionless contact interfaces, \href{https://doi.org/10.1016/j.cma.2023.116641}{Comput. Methods Appl. Mech. Eng.}, 419 (2024) 116641.
	
	\bibitem {JW50}J. Zhang, H. Dai, Q. Wang, Highly accurate reconstruction harmonic balance scheme for Three-Body libration point orbits, \href{https://doi.org/10.2514/1.G008551}{J. Guid. Control Dyn.}, 48 (2025) 1148-1157.
	
	\bibitem {JW51}J. Jiao, J. Xu, X. Yuan, L. Chen, Axisymmetric 3 : 1 internal resonance of thin-walled hyperelastic cylindrical shells under both axial and radial excitations, \href{https://doi.org/10.1007/s10409-022-09006-x}{Acta Mech. Sin.}, 38 (2022) 521417.
	
	\bibitem {JW36}J. Kappauf, S. Bäuerle, H. Hetzler, A combined FD-HB approximation method for steady-state vibrations in large dynamical systems with localised nonlinearities, \href{https://doi.org/10.1007/s00466-022-02225-3}{Comput. Mech.}, 70 (2022) 1241-1256.
	
	\bibitem {JW37}S. Karkar, B. Cochelin, C. Vergez, A comparative study of the harmonic balance method and the orthogonal collocation method on stiff nonlinear systems, \href{https://doi.org/10.1016/j.jsv.2014.01.019}{J. Sound Vib.}, 333 (2014) 2554-2567.
	
	\bibitem {JW38}L.C. Young, Orthogonal collocation revisited, \href{https://doi.org/10.1016/j.cma.2018.10.019}{Comput. Methods Appl. Mech. Eng.}, 345 (2019) 1033-1076.
	
	\bibitem {JW39}J. Wu, L. Hong, Y. Xu, J. Jiang, An efficient and robust approach for continuation and bifurcation analysis of quasi-periodic solutions by multi-harmonic balance method, \href{https://doi.org/10.1016/j.jsv.2025.118943}{J. Sound Vib.}, (2025) 118943.
	
	\bibitem {JW40}M. Li, H. Yan, L. Wang, Nonlinear model reduction for a cantilevered pipe conveying fluid: A system with asymmetric damping and stiffness matrices, \href{https://doi.org/10.1016/j.ymssp.2022.109993}{Mech. Syst. Sig. Process.}, 188 (2023) 109993.
	
	\bibitem {JW52}M. Li, G. Haller,Nonlinear analysis of forced mechanical systems with internal resonance using spectral submanifolds, Part II: Bifurcation and quasi-periodic response, \href{https://doi.org/10.1007/s11071-022-07476-6}{Nonlinear Dyn.}, 110 (2022) 1045-1080.
	
	\bibitem {JW41}S.H. Strogatz, Nonlinear dynamics and chaos: with applications to physics, biology, chemistry, and engineering, (CRC, Boca Raton, 2015), \href{https://doi.org/10.1201/9780429492563}{https://doi.org/10.1201/9780429492563}.
	
	\bibitem {JW42}Y. Colaïtis, A. Batailly, Stability analysis of periodic solutions computed for blade-tip/casing contact problems, \href{https://doi.org/10.1016/j.jsv.2022.117219}{J. Sound Vib.}, 538 (2022) 117219.
	
	\bibitem {JW43}J. Wu, L. Hong, J. Jiang, A robust and efficient stability analysis of periodic solutions based on harmonic balance method and Floquet-Hill formulation, \href{https://doi.org/10.1016/j.ymssp.2022.109057}{Mech. Syst. Sig. Process.}, 173 (2022) 109057.
	
	\bibitem {JW44}L. Guillot, A. Lazarus, O. Thomas, C. Vergez, and B. Cochelin, A purely frequency based Floquet-Hill formulation for the efficient stability computation of periodic solutions of ordinary differential systems, \href{https://doi.org/10.1016/j.jcp.2020.109477}{J. Comput. Phys.}, 416 (2020) 109477.
	
	\bibitem {JW45}F. Bayer, R. Leine, Sorting-free Hill-based stability analysis of periodic solutions through Koopman analysis, \href{https://doi.org/10.21203/rs.3.rs-2183060/v1}{Nonlinear Dyn.}, 111 (2023) 8439-8466.
	
	\bibitem {JW46}J.H. Argyris, G. Faust, M. Haase, and R. Friedrich, An exploration of dynamical systems and chaos, (Springer, Heidelberg, 2015), \href{https://doi.org/10.1007/978-3-662-46042-9}{https://doi.org/10.1007/978-3-662-46042-9}.
	
	\bibitem {JW47}G. Benettin, L. Galgani, A. Giorgilli, and J.-M. Strelcyn, Lyapunov Characteristic Exponents for smooth dynamical systems and for hamiltonian systems; A method for computing all of them. Part 2: Numerical application, \href{https://doi.org/10.1007/BF02128237}{Meccanica}, 15 (1980) 21-30.
	
	\bibitem {JW48}R. Fiedler, H. Hetzler, S. Bäuerle, Efficient numerical calculation of Lyapunov-exponents and stability assessment for quasi-periodic motions in nonlinear systems, \href{https://doi.org/10.1007/s11071-024-09497-9}{Nonlinear Dyn.}, 112 (2024) 8299-8327.
	
	
	
	
\end{thebibliography}
\end{document}